\documentclass[oneside,english,reqno]{amsart}
\usepackage{latexsym,amsmath,amssymb,amsfonts,amsthm,verbatim}
\usepackage[utf8]{inputenc}
\usepackage{color}
\usepackage{eufrak}
\usepackage{cases}

\makeatletter

  \theoremstyle{plain}
  \newtheorem*{thm*}{Theorem}
  \newtheorem{corollary}{Corollary}
  \theoremstyle{plain}
  \newtheorem{lemma}{Lemma}
  \theoremstyle{remark}
  \newtheorem{remark}{Remark}
 \theoremstyle{definition}
 \newtheorem*{defn*}{Definition}

\newcommand{\eqdef}{\overset{{\mbox{\tiny  def}}}{=}}
\def\hp{\mathbf{H}}
\def\Ga{\Gamma}
\def\vp{\varphi}

\makeatother
\begin{document}

\title{Malliavin Calculus and Self Normalized Sums\\ \today}

\author{Solesne Bourguin and Ciprian A. Tudor}
\address{\noindent SAMM, EA 4543,
 Universit\'e Paris 1 Panth\'eon Sorbonne,
90 Rue de Tolbiac, 75634 Paris Cedex France}
\address{Laboratoire Paul Painlev\'e, Universit\'e de Lille 1,
 F-59655 Villeneuve d'Ascq, France}

\email{solesne.bourguin@univ-paris1.fr and tudor@math.univ-lille1.fr }

\begin{abstract}
We study the self-normalized sums of independent random variables from the perspective of the Malliavin calculus. We give the chaotic expansion for them and we prove a Berry-Ess\'een bound  with respect to several distances.
\end{abstract}
\maketitle

\vskip0.3cm

{\bf 2010 AMS Classification Numbers:}  60F05, 60H05, 91G70.

 \vskip0.3cm


\section{Introduction}

\noindent Let $(\Omega, {\mathcal{F}}, P)$ be a probability space and $(W_{t})_{t \geq 0}$ a Brownian motion on this space. Let $F$ be a random variable defined on $\Omega$ which is differentiable in the sense of the Malliavin calculus.
Then using  the so-called Stein's method introduced by Nourdin and Peccati in \cite{NoPe1} (see also \cite{NoPe2} and \cite{NoPe3}), it is possible to measure the distance between the law of $F$ and the standard normal law $\mathcal{N}(0,1)$. This distance can be defined in several ways, such as the Kolmogorov distance, the Wasserstein distance, the total variation distance or the Fortet-Mourier distance. More precisely we have, if ${\mathcal{L}}(F)$  denotes the law of $F$,
\begin{equation*}
d({\mathcal{L}}(F), \mathcal{N}(0,1)) \leq c\sqrt{ \mathbf{E}\left( 1-\langle DF , D(-L) ^{-1} F \rangle _{L^{2}([0,1])}\right)^{2}}.
\end{equation*}
Here $D$ denotes the Malliavin derivative with respect to $W$, and $L$ is the generator of the Ornstein-Uhlenbeck semigroup. We will explain in the next section  how these operators are defined. The constant $c$ is equal to 1 in the case of the Kolmogorov distance as well as in the case of the Wasserstein distance, $c=2$ for the  total variation distance and $c=4$ in the case of the Fortet-Mourier distance.

\noindent Our purpose is to apply these techniques to self-normalized sums. Let us recall some basic facts on this topic. We refer to \cite{DLS} and the references therein for a more detailed exposition. Let $X_{1},X_{2}, \cdots $ be independent random variables. Set $S_{n}=\sum_{i=1}^{n}X_{i} $ and $V_{n}^{2} =\sum_{i=1}^{n} X_{i}^{2}$. Then $\frac{S_{n}}{V_{n}}$ converges in distribution as $n\to \infty$ to the standard normal law $\mathcal{N}(0,1)$ if and only if $\mathbf{E}(X)=0$ and $X$ is in the domain of attraction of the standard normal law (see \cite{DLS}, Theorem 4.1). The ``if'' part of the theorem has been known for a long time (it appears in \cite{Ma}) while the ``only if'' part remained open until its proof in \cite{Gi}.
The Berry-Ess\'een theorem for self-normalized sums has been also widely studied. We refer to \cite{BK2} and \cite{Shao} (see also \cite{BK1}, \cite{BK3} for the situation where the random variables $X_{i}$ are non i.i.d. ).  These results say that  the Kolmogorov distance between the law of $\frac{S_{n}}{V_{n}}$ and the standard normal law is less than
\begin{equation*}
C\left( B_{n}^{-2}\sum_{i=1}^{N} \mathbf{E}\left( X_{i}^{2} 1_{(\vert X_{i}\vert >B_{n} )}\right) + B_{n}^{-3}\sum_{i=1}^{N} \mathbf{E}\left( X_{i}^{3} 1_{(\vert X_{i}\vert \geq B_{n} )}\right) \right)
 \end{equation*}
 where $B_{n}=\sum_{i=1}^{n} \mathbf{E}(X_{i}^{2})$ and $C$ is an absolute constant.
 We mention that, as far as we know, these results only exist for the Kolmogorov distance.
 To use our techniques based on the Malliavin calculus and multiple stochastic integrals, we will put ourselves on a Gaussian space where we will consider the following particular case: the random variables $X_{i}$ are the increments of the Wiener process $X_{i}=W_{i}-W_{i-1}$. The Berry-Ess\'een bound from above reduces to (see \cite{DLS}, page 53): for $2<p\leq 3$
\begin{equation}
\label{borneDLS}
\sup_{z\in \mathbb{R}} \left| \mathbf{P}(F_{n} \leq z) -\Phi (z) \right| \leq 25 \mathbf{E}\left( |Z| ^{p}\right)  n^{1-\frac{p}{2} }
\end{equation}
where $Z$ is a standard normal random variable and $\Phi$ is its repartition function. In particular for $p=3$ we get
\begin{equation}\label{borneDLS1}
\sup_{z\in \mathbb{R}} \left| \mathbf{P}(F_{n} \leq z) -\Phi (z) \right| \leq 25 \mathbf{E}\left(|Z| ^{3}\right)  n^{-\frac{1}{2} }.
\end{equation}
 We will compare our result with the above relation (\ref{borneDLS1}). The basic idea is as follows: we are able to find the chaos expansion into multiple Wiener-It\^o integrals  of the random variable $\frac{S_{n}}{V_{n}}$ for every $n\geq 2$ and to compute its Malliavin derivative. Note that the random variable $\frac{S_{n}}{V_{n}}$ has a decomposition into an infinite sum of multiple integrals in contrast to the examples provided in the papers \cite{BoTu}, \cite{NoPe1}, \cite{NoPe2}.  Then we compute the Berry-Ess\'een bound given by $\sqrt{ \mathbf{E}\left( 1-\langle DF , D(-L) ^{-1} F \rangle _{L^{2}([0,1])}\right)^{2}}$ by using properties of multiple stochastic integrals.  Of course, we cannot expect to obtain a rate of convergence better than $c\frac{1}{\sqrt{n}}$, but we have an explicit (although complicated) expression of the constant appearing in this bound and our method is available for several distances between the laws of random variables (not limited to the Kolmogorov distance). This aspect of the problem seems to be new. This computation of the Berry-Ess\'een bound is also interesting in and of itself as it brings to light original relations involving Gaussian measure and Hermite polynomials. It gives an exact expression of the chaos expansion of the self normalized sum and  it also shows that the convergence to the normal law of $\frac{S_{n}}{V_{n}}$ is uniform with respect to the chaos, in the sense that every chaos of  $\frac{S_{n}}{V_{n}}$ is convergent to the standard normal law and that the rate is the same for every chaos.

\noindent We have organized our paper as follows: Section 2 contains the elements of the Malliavin calculus needed in the paper and in Section 3 we discuss the chaos decomposition of self-normalized sums as well as study the asymptotic behavior of the coefficients appearing in this expansion. Section 4 contains the computation of the Berry-Ess\'een bound given in terms of the Malliavin calculus.

\section{Preliminaries}

\noindent We will begin by describing the basic tools of multiple Wiener-It\^o integrals and Malliavin calculus that will be needed in our paper. Let $(W_{t})_{t\in [0,T] }$ be a classical Wiener process on a standard Wiener space $\left( \Omega , {\mathcal{F}}, P\right)$. If $f\in L^{2}([0,T]^{n}) $ with $n\geq 1$ integer, we introduce the multiple Wiener-It\^o integral of $f$ with respect to $W$. We refer to \cite{N} for a detailed exposition of the construction and the properties of multiple Wiener-It\^o integrals.

\noindent Let $f\in {\mathcal{S}}_{n}$, which means that there exists $n \geq 1$ integers such that
\begin{equation*}
f:=\sum_{i_{1}, \cdots  , i_{n}}c_{i_{1}, \cdots  ,i_{n} }1_{A_{i_{1}}\times \cdots  \times A_{i_{n}}}
\end{equation*}
where the coefficients satisfy $c_{i_{1}, \cdots  ,i_{n} }=0$ if two indices $i_{k}$ and $i_{\ell}$ are equal and the sets
 $A_{i}\in {\mathcal{B}}([0,T]) $ are disjoints. For a such step function $f$ we define

\begin{equation*}
I_{n}(f):=\sum_{i_{1}, \cdots  , i_{n}}c_{i_{1}, \cdots  ,i_{n} }W(A_{i_{1}})\cdots  W(A_{i_{n}})
\end{equation*}
where we put $W([a,b])= W_{b}-W_{a}$. It can be seen that the application $I_{n}$ constructed above from ${
\mathcal{S}}_{n}$ equipped with the scaled norm $\frac {1}{\sqrt{n!}}\Vert\cdot\Vert_{L^{2}([0, T] ^{n})}$ to $L^{2}(\Omega)$ is
an isometry on ${\mathcal{S}}_{n}$, i.e. for $m,n$ positive integers,
\begin{eqnarray*}
\mathbf{E}\left(I_{n}(f) I_{m}(g) \right) &=& n! \langle f,g\rangle _{L^{2} ([0,T]^{n})}\quad \mbox{if } m=n,\\
\mathbf{E}\left(I_{n}(f) I_{m}(g) \right) &= & 0\quad \mbox{if } m\not=n.
\end{eqnarray*}
It also holds that
\begin{equation*}
I_{n}(f) = I_{n}\big( \tilde{f}\big)
\end{equation*}
where $\tilde{f} $ denotes the symmetrization of $f$ defined by $$\tilde{f}%
(x_{1}, \cdots  , x_{x}) =\frac{1}{n!} \sum_{\sigma \in {\mathcal S}_{n}} f(x_{\sigma (1) }, \cdots  , x_{\sigma (n) } ). $$ Since the set ${\mathcal{S}}_{n}$ is dense in $L^{2}([0,T]^{n})$ for every $n\geq  2$, the mapping $I_{n}$ can be extended to an isometry from $L^{2}([0,T]^{n})$ to  $L^{2}(\Omega)$ and the above properties hold true for this extension. Note also that $I_{n}$ can be viewed as an iterated stochastic integral (this follows e.g. by It\^o's formula)
\begin{equation*}
I_{n}(f)=n!\int_{0}^{1}\int_{0} ^{t_{n}} \cdots  \int_{0}^{t_{2}}f(t_{1}, \cdots  , t_{n}) dW_{t_{1}}\cdots  dW_{t_{n}}
\end{equation*}

\noindent We recall the product for two multiple integrals (see \cite{N}): if $f\in L^{2}([0,T] ^{n})$ and $g\in L^{2}([0,T]  ^{m} )$ are symmetric, then it holds that
\begin{equation}  \label{prod}
I_{n}(f)I_{m}(g)= \sum_{\ell=0}^{m\wedge n} \ell! C_{m}^{\ell}C_{n}^{\ell} I_{m+n-2\ell} (f\otimes _{\ell}g)
\end{equation}
where the contraction $f\otimes _{\ell}g$ belongs to $L^{2}([0,T]^{m+n-2\ell}) $ for $\ell=0, 1, \cdots  , m\wedge n$ and is
given by
\begin{eqnarray*}\nonumber
&&(f\otimes _{\ell} g) ( s_{1}, \cdots  , s_{n-\ell}, t_{1}, \cdots  , t_{m-\ell})  \\
&=&\int_{[0,T] ^{\ell} } f( s_{1}, \cdots  , s_{n-\ell}, u_{1}, \cdots  ,u_{\ell})g(t_{1}, \cdots  , t_{m-\ell},u_{1},
\cdots  ,u_{\ell}) du_{1}\cdots  du_{\ell}.
\end{eqnarray*}
We recall that any square integrable random variable that is measurable with respect to the $\sigma$-algebra generated by $W$ can be expanded into an orthogonal sum of multiple stochastic integrals
\begin{equation}\label{fin1}
F=\sum_{n\geq 0}I_{n}(f_{n})
\end{equation}
where $f_{n}\in L^{2}([0,1]^{n})$ are (uniquely determined) symmetric functions and $I_{0}(f_{0})=\mathbf{E}\left( F\right)$.

\noindent Let $L$ be the Ornstein-Uhlenbeck operator
\begin{equation*}
LF=-\sum_{n\geq 0} nI_{n}(f_{n}) \mbox{ and } L^{-1} F = -\sum_{n\geq 1} \frac{1}{n}  I_{n}(f_{n})
\end{equation*}
if $F$ is given by (\ref{fin1}). We denote by $D$  the Malliavin  derivative operator that acts on smooth functionals of the form $F=g(W(\varphi _{1}), \cdots  , W(\varphi_{n}))$ where $g$ is a smooth function with compact support and $\varphi_{i} \in L^{2}([0,1])$. For $i=1,\cdots ,n$, the derivative operator is defined by
\begin{equation*}
DF=\sum_{i=1}^{n}\frac{\partial g}{\partial x_{i}}(B(\varphi _{1}), \cdots  , B(\varphi_{n}))\varphi_{i}.
\end{equation*}
The operator $D$ can be extended to the closure $\mathbb{D}^{p,2}$ of smooth functionals with respect to the norm
\begin{equation*}
\Vert F\Vert _{p,2}^{2} = \mathbf{E}\left( F^{2}\right) + \sum_{i=1}^{p} \mathbf{E}\left(  \Vert D^{i} F\Vert ^{2} _{L^{2}([0,1]^{i})}\right)
\end{equation*}
where the $i^{\mbox{\tiny{th}}}$ order Malliavin derivative $D^{i}$ is defined iteratively.

\noindent Let us recall how this  derivative acts  for random variables in a finite chaos.
If $f\in L^{2}([0,T]^{n}) $  is a symmetric function, we will use the following rule to differentiate in the Malliavin sense
\[
D_{t}I_{n}(f)=n \, I_{n-1}(f(\cdot,t)),\hskip0.5cmt\in \mathbb{R}.
\]

\noindent Let us also recall how the distances between the laws of random variables are defined. We have
\begin{equation*}
d({\mathcal{L}}(X),{\mathcal{L}}(Y)) =\sup_{h\in {\mathcal{A}}}\left(  \left| \mathbf{E}\left( h(X)\right) - \mathbf{E}\left( h(Y)\right) \right| \right)
\end{equation*}
where ${\mathcal{A}}$ denotes a set of functions. When ${\mathcal{A}}=\{ h: \Vert h\Vert _{L}\geq 1\}$ (here $\Vert \cdot \Vert _{L}$ is the Lipschitz norm) we obtain the Wasserstein distance,  when ${\mathcal{A}}=\{ h: \Vert h\Vert _{BL}\geq 1\}$ (with $\Vert \cdot \Vert _{LB} = \Vert \cdot \Vert _{L} + \Vert \cdot \Vert _{\infty}$) we get the Fortet-Mourier distance, when ${\mathcal{A}}$ is the set of indicator functions of Borel sets we obtain the total variation distance, and when ${\mathcal{A}}$ is the  set  of indicator functions of the form $1_{(-\infty, z)}$ with $z\in \mathbb{R}$, we obtain the Kolmogorov distance that has been presented  above.
\vskip0.3cm

\section{Chaos decomposition of self-normalized sums}

\noindent The tools of the Malliavin calculus presented above can be successfully applied in order to study self-normalized sums. Because of the nature of Malliavin calculus, we put ourselves in a Gaussian setting and we consider $X_{i}= W_{i}-W_{i-1}$ to be the increments of a classical Wiener process $W$. We then consider the sums
\begin{equation*}
S_{n} =\sum_{i=1}^{n} X_{i} \hskip0.5cm \mbox{and}\hskip0.5cm  V_{n}^{2}= \sum_{i=1}^{n} X_{i}^{2}
\end{equation*}
as well as the \textit{self-normalized sum} $F_{n}$ defined by
\begin{equation}\label{fn}
F_{n}= \frac{S_{n}}{V_{n}}= \frac{W_{n}}{\left( \sum_{i=1}^{n}(W_{i+1} -W_{i})^{2} \right) ^{\frac{1}{2} }}.
\end{equation}

\noindent Let us now concentrate our efforts on finding the chaotic decomposition of the random variable $F_{n}$. This will be the key to computing Berry-Ess\'een bounds for the distance between the law of $F_{n}$ and the standard normal law in the next section.
\begin{lemma}
\label{lemmeChaosExp}
Let $F_{n}$ be given by (\ref{fn}) and let $f: \mathbb{R}^{n}\to \mathbb{R}$ be given by
\begin{equation}\label{ff}
f(x_{1}, \cdots , x_{n}) = \frac{x_{1}+\cdots  + x_{n}}{(x_{1}^{2}+ \cdots  + x_{n}^{2} ) ^{\frac{1}{2}}}.
\end{equation}
Then for every $n\geq 2$, we have
\begin{equation*}
F_{n}=\sum_{k\geq 0} \frac{1}{k!}\sum_{i_{1}, \cdots , i_{k}=1}^{n} a_{i_{1},\cdots , i_{k}} I_{k} \left( \varphi_{i_{1}}\otimes  \cdots \otimes \varphi_{i_{k}}\right)
\end{equation*}
with
\begin{equation}\label{eqn:ai}
a_{i_{1},\cdots , i_{k}} \eqdef \mathbf{E} \left( \frac{\partial ^{k} f}{\partial x_{i_{1}},\cdots ,x_{i_{k}}}\left( W(\varphi _{1}), \cdots , W(\varphi_{n}) \right)\right).
\end{equation}
\end{lemma}

\noindent {\bf Proof: }We use the so-called Stroock's formula (see \cite{N}). The Wiener chaos expansion of a smooth  (in the sense of Malliavin calculus) random variable $F$ is  given by
\begin{equation}\label{eqn:formStr}
F=\sum_{k\geq 0} \frac{1}{k!} I_{k} \left( \mathbf{E} \left( D^{k}F\right) \right)
\end{equation}
where $D^{k}$ denotes the $k^{\mbox{\tiny{th}}}$  iterated Malliavin derivative.
Note that $F_{n}$ can be written as
\begin{equation*}
F_{n} = f\left( W(\varphi _{1}), \cdots , W(\varphi_{n}) \right)
\end{equation*}
where
\begin{equation*}
\varphi_{i} = 1_{[i-1,i]}, \hskip0.3cm i=1,\cdots ,n.
\end{equation*}

\noindent The chain rule for the  Malliavin derivative  (with  $D_{s}W(\varphi)= \varphi (s)$) yields
\begin{equation*}
DF_{n}=\sum_{i=1}^{n} \frac{\partial f}{\partial x_{i}} \left( W(\varphi _{1}), \cdots , W(\varphi_{n}) \right) \varphi_{i}
\end{equation*}
and proceeding recursively leads to the formula
\begin{equation*}
D^{k}F_{n}= \sum_{i_{1}, \cdots , i_{k}=1}^{n} \frac{\partial ^{k} f}{\partial x_{i_{1}},\cdots ,x_{i_{k}}}\left( W(\varphi _{1}), \cdots , W(\varphi_{n}) \right)\varphi_{i_{1}}\otimes  \cdots \otimes \varphi_{i_{k}}.
\end{equation*}
Thus we obtain
\begin{eqnarray*}
I_{k}\big(\mathbf{E} (D^{k}F_{n}) \big) &=&  \sum_{i_{1}, \cdots , i_{k}=1}^{n} \mathbf{E}\,\left( \frac{\partial ^{k} f}{\partial x_{i_{1}},\cdots ,x_{i_{k}}}\left( W(\varphi _{1}), \cdots , W(\varphi_{n}) \right)\right) I_{k} \left( \varphi_{i_{1}}\otimes  \cdots \otimes \varphi_{i_{k}}\right)\\
&=& \sum_{i_{1}, \cdots , i_{k}=1}^{n} a_{i_{1},\cdots , i_{k}} I_{k} \left( \varphi_{i_{1}}\otimes  \cdots \otimes \varphi_{i_{k}}\right)
\end{eqnarray*}
where $a_{i_{1},\cdots , i_{k}}$ are defined by
(\ref{eqn:ai}). Thus  from \eqref{eqn:formStr} it follows that,
\begin{equation*}
F_{n}=\sum_{k\geq 0} \frac{1}{k!}\sum_{i_{1}, \cdots , i_{k}=1}^{n} a_{i_{1},\cdots , i_{k}} I_{k} \left( \varphi_{i_{1}}\otimes  \cdots \otimes \varphi_{i_{k}}\right).
\end{equation*}
\qed

\begin{remark}
The coefficients $a_{i_{1},\cdots , i_{k}}$ also depend on $n$. We omit $n$ in their notation in order to simplify the presentation.
\end{remark}

\subsection{Computing the coefficients in the chaos expansion}
In this subsection, we explicitly compute the coefficients $a_{i_{1},\cdots  , i_{k}}$ appearing in Lemma \ref{lemmeChaosExp}. Let $\hp_n(x)$ denote the $n^{\text{th}}$ Hermite polynomial:
\begin{equation*}
\hp_n(x) = (-1)^n e^{x^2/2} {d^n \over dx^n} e^{-x^2/2}.
\end{equation*}
Define
\begin{eqnarray*}
W_n &\eqdef & W (\varphi_1) + W(\varphi_2)  + \cdots  + W(\varphi_n) \; \\
V_n &\eqdef &  \Big(\sum_{i=1}^n W (\varphi_i)^2 \Big)^{1/2}
\end{eqnarray*}
Let us first give the following lemma that can be proved using integration by parts.
\begin{lemma}\label{lem:intbypart}
 For every $1\leq i_{1},., i_{k} \leq n$, let $a_{i_{1}, \cdots  i_{k}}$
 be as defined in \eqref{eqn:ai}. Let $d_r, 1 \leq r \leq n$ denote
 the number of times the integer $r$ appears in the sequence $\{i_1,i_2,\cdots ,i_k\}$ with
 $\sum_{r=1}^n d_r = k$. Then we have
\begin{equation*}
a_{i_{1}, \cdots  i_{k}}= \mathbf{E} \left( {W_{n} \over V_n} \prod_{r=1}^n \hp_{d_r}\big(W(\varphi _{r})\big)  \right).
\end{equation*}
\end{lemma}
\noindent {\bf Proof: }If $X \sim \mathcal{N}(0,1)$, then for any $g \in C^{(n)}(\mathbb{R})$ with $g$ and its derivatives having polynomial growth at infinity, we have the Gaussian integration by parts formula
\begin{equation*}
\mathbf{E}(g^{(n)}(X)) = \mathbf{E}(g(X) \hp_n(X)).
\end{equation*}
where $g^{(n)}(x) \eqdef {d^n \over dx^n} g(x)$.

\noindent Notice that the function $f$ defined in (\ref{ff})
satisfies $|f(x)| \leq C|x|, \forall x \in \mathbb{R}^n$ for a constant $C$, and
thus applying the above integration by parts formula recursively yields
\begin{eqnarray*}
a_{i_{1}, \cdots  i_{k}}&=&\frac{1}{(\sqrt{2\pi}) ^{n}}\int_{\mathbb{R}^{n}} \left( \frac{\partial ^{k} f}{\partial x^{d_1}_{1},\cdots ,x^{d_n}_{n}}\right) (x_{1},\cdots , x_{n}) \, e^{-\frac{x_{1} ^{2}}{2}} \cdots
e^{-\frac{x_{n} ^{2}}{2} }dx_{1}\cdots dx_{n}\\
&=&\frac{1}{(\sqrt{2\pi}) ^{n}}\int_{\mathbb{R}^{n}} \left( \frac{\partial ^{k} f}{\partial x^{d_1}_{1},\cdots ,x^{d_{n-1}}_{n-1}}\right) (x_{1},\cdots , x_{n})\,\,\hp_{d_n}(x_n) \,e^{-\frac{x_{1} ^{2}}{2}} \cdots
e^{-\frac{x_{n} ^{2}}{2} }\,dx_{1}\cdots dx_{n}\\
&=& \frac{1}{(\sqrt{2\pi}) ^{n}}\int_{\mathbb{R}^{n}} f(x_{1},\cdots , x_{n}) \prod_{r=1}^n \hp_{d_r}(x_r)\,e^{-\frac{x_{1} ^{2}}{2}} \cdots
e^{-\frac{x_{n} ^{2}}{2} }\,dx_{1}\cdots dx_{n}\\
&=&\mathbf{E} \left( {W_{n} \over V_n} \prod_{r=1}^n \hp_{d_r}\big(W(\varphi _{r})\big)  \right).
\end{eqnarray*}
This concludes the proof of the Lemma. \qed \\

\noindent The next step in the calculation of the coefficient is to notice that $a_{i_{1}, \cdots  i_{k}} = 0$ when $k$ is even. This is the object of the following Lemma.
\begin{lemma} \label{lem:evenzero} If $k$ is even, then
\begin{equation*}
a_{i_{1}, \cdots  i_{k}}=0.
\end{equation*}
\end{lemma}
\noindent {\bf Proof: }Let $k$ be an even number and $d_1, d_2, \cdots , d_n$ be as defined in Lemma \ref{lem:intbypart}.  By Lemma \ref{lem:intbypart}, we have
\begin{equation}
\label{eqn:aievenzerotrick}
a_{i_{1}, \cdots  i_{k}} = \sum_{u=1}^{n} \mathbf{E} \left( \frac{W(\varphi_u)}{V_{n}} \prod_{r=1}^n \hp_{d_r}\big(W(\varphi _{r})\big)  \right).
\end{equation}
Note that the product $\prod_{r=1}^n \hp_{d_r}\big(W(\varphi _{r}))$ is an even function of
$(W(\varphi_1), W(\varphi_2), \cdots , W(\varphi_n))$. Indeed, since $k$ is even and
$\sum_{r=1}^n d_r = k$, either all of the integers $d_r, r\leq n$ are even or there is an even
number of odd integers in $d_r, r\leq n$. In either case the product $\prod_{r=1}^n \hp_{d_r}\big(W(\varphi _{r}))$
is an even function of $(W(\varphi_1), W(\varphi_2), \cdots , W(\varphi_n))$, since $\hp_m(x) = \hp_m(-x)$ for all even $m \in \mathbb{N}$ and $\hp_m(x) = -\hp_m(-x)$ for all odd $m \in \mathbb{N}$.

\noindent Thus for each $u \leq n$, the expression $ \frac{W(\varphi_u)}{V_{n}} \prod_{r=1}^n \hp_{d_r}\big(W(\varphi _{r})\big)$
is an odd function of $W(\varphi_u)$ and thus has expectation zero since $W(\phi_u)$
is a standard Gaussian random variable. The fact that \eqref{eqn:aievenzerotrick} is a sum of such expectations concludes the proof.
\qed \vskip0.3cm

\noindent As a consequence of Lemma \ref{lem:evenzero}, we have
\begin{equation}\label{eqn:chaosFN1}
F_{n}=\sum_{k\geq 0} \frac{1}{(2k+1)!}\sum_{i_{1}, \cdots , i_{2k+1}=1}^{n} a_{i_{1},\cdots , i_{2k+1}} I_{2k+1} \left( \varphi_{i_{1}}\otimes  \cdots \otimes \varphi_{i_{2k+1}}\right).
\end{equation}
This implies that in order to compute the coefficients $a_{i_{1}, \cdots  i_{k}}$, it suffices to focus on the case where $k$ is odd. Before stating the first result in this direction, let us give the following technical lemma.
\begin{lemma}
\label{lemmaxx}
Let $k \geq 0$ be a positive integer and let $d_r, 1 \leq r \leq n$ denote the number of times the integer $r$ appears in the sequence $\{i_1,i_2,\cdots ,i_{2k+1}\}$ with $\sum_{r=1}^n d_r = 2k+1$. Then, if there is more than one odd integer in the sequence $d_r, 1 \leq r \leq n$, for each $1 \leq i \leq n$,
\begin{equation*}
\mathbf{E}\left[ \frac{1}{V_{n}}W(\varphi_{i} ) \hp_{d_{1}} \left(  W(\varphi_{1})\right) \hp_{d_{2}} \left(  W(\varphi_{2})\right)\cdots  \hp_{d_{n} } \left(  W(\varphi_{n})\right)\right] = 0.
\end{equation*}
\end{lemma}
\noindent {\bf Proof: }Note that the equality $\sum_{r=1}^n d_r = 2k+1$ implies that there can only be an odd number of odd integers in the sequence $d_r$, otherwise the sum $\sum_{r=1}^n d_r$ could not be odd. Therefore, more than one odd integer in the sequence $d_r$ means that there are at least three of them. We will prove the Lemma for this particular case of three odd integers in the sequence $d_r$ for the sake of readability of the proof, as the other cases follow with the exact same arguments. Hence, assume that there are three odd integers $d_i$, $d_k$ and $d_l$ in the sequence $d_r, 1 \leq r \leq n$. We will first consider the case where $i$ is different than $j,k,l$. Then,
\begin{eqnarray*}
&& \mathbf{E}\left[ \frac{1}{V_{n}}W(\varphi_{i} ) \hp_{d_{1}} \left(  W(\varphi_{1})\right) \hp_{d_{2}} \left(  W(\varphi_{2})\right)\cdots  \hp_{d_{n} } \left(  W(\varphi_{n})\right)\right] \\
&&= \frac{1}{(2n)^{\frac{n}{2}}}\int_{\mathbb{R}^{n}}\frac{x_{i} \hp_{d_{1}} \left(x_1\right) \cdots  \hp_{d_{n}} \left(x_n\right)}{\sqrt{x_{1}^{2} + \cdots  + x_{n}^{2}}}e^{-\frac{1}{2}(x_{1}^{2} + \cdots  + x_{n}^{2})}dx_{1} \cdots  dx_{n} \\
&& =\frac{1}{(2n)^{\frac{n}{2}}}\int_{\mathbb{R}^{n-1}}x_{i} \hp_{d_{1}} \left(x_1\right) \cdots  \hp_{d_{j-1}} \left(x_{j-1}\right)\hp_{d_{j+1}} \left(x_{j+1}\right) \cdots  \hp_{d_{n}} \left(x_n\right) \\
&& \times \left(\int_{\mathbb{R}}\frac{\hp_{d_{j}} \left(x_j\right)}{\sqrt{x_{1}^{2} + \cdots  + x_{n}^{2}}}e^{-\frac{x_{j}^{2}}{2}}dx_j \right) \mbox{exp}\left[ -\frac{1}{2}\sum_{\underset{p \neq j }{p=1}}^{n}x_{p}^{2}\right] dx_{1} \cdots  dx_{j-1}dx_{j+1} \cdots  dx_{n}.
\end{eqnarray*}
$d_j$ beeing odd, $\hp_{d_j}$ is an odd function of $x_j$ and $x_j \mapsto \frac{\hp_{d_{j}} \left(x_j\right)}{\sqrt{x_{1}^{2} + \cdots  + x_{n}^{2}}}e^{-\frac{x_{j}^{2}}{2}}$ is also an odd function of $x_j$. Thus, $\int_{\mathbb{R}}\frac{\hp_{d_{j}} \left(x_j\right)}{\sqrt{x_{1}^{2} + \cdots  + x_{n}^{2}}}e^{-\frac{x_{j}^{2}}{2}}dx_j = 0$ and finally
\begin{equation*}
\mathbf{E}\left[ \frac{1}{V_{n}}W(\varphi_{i} ) \hp_{d_{1}} \left(  W(\varphi_{1})\right) \hp_{d_{2}} \left(  W(\varphi_{2})\right)\cdots  \hp_{d_{n} } \left(  W(\varphi_{n})\right)\right] = 0.
\end{equation*}
The other cases one could encounter is when $i = j$ or $i = k$ or $i = l$ and the proof follows based on the exact same argument.
\qed \vskip0.3cm

\noindent We can now state the following key result that will allow us to perform further calculations in order to explicitly determine the coefficients $a_{i_{1}, \cdots  i_{k}}$.

\begin{lemma}
For every $k \geq 0$ and for every $1 \leq i_1,\cdots ,i_{2k+1} \leq n$, let $d_{r}^{\star}, 1 \leq r \leq n$ be the number of times the integer $r$ appears in the sequence $\{i_1,\cdots ,i_{2k+1}\}$. Then,
\begin{equation}
\label{aoddtemp}
a_{1_{1},\cdots , i_{2k+1}}= \mathbf{E}\left[ \frac{1}{V_{n}}W(\varphi_{1} ) \hp_{d_{1}^{\star}} \left(  W(\varphi_{1})\right) \hp_{d_{2}^{\star}} \left(  W(\varphi_{2})\right)\cdots  \hp_{d_{n}^{\star}} \left(  W(\varphi_{n})\right)\right]
\end{equation}
if there is only one odd integer in the sequence $d_{r}^{\star}, 1 \leq r \leq n$. If there is more than one odd integer in the sequence $d_{r}^{\star}, 1 \leq r \leq n$, we have $a_{1_{1},\cdots , i_{2k+1}}= 0$.
\end{lemma}
\begin{remark}
Note that in (\ref{aoddtemp}), it might be understood that $d_1^{\star}$ is always the only odd integer in $d_{r}^{\star}, 1 \leq r \leq n$. This is obviously not always the case and if $d_1^{\star}$ is not the odd integer but let's say, $d_i^{\star}$ with $1 < i \leq n$ is, one can use the equality in law between $W(\varphi_{i} )$ and $W(\varphi_{1} )$ to perform an index swap $(i\leftrightarrow 1)$ and the equality (\ref{aoddtemp}) remains unchanged.
\end{remark}
\begin{remark}
If one is in the case where $a_{1_{1},\cdots , i_{2k+1}} \neq 0$, one can rewrite $d_1^{\star},d_2^{\star},\cdots ,d_n^{\star}$ as $2d_{1}+1,2d_2,\cdots ,2d_n$ and finally rewrite (\ref{aoddtemp}) as
\begin{equation}
\label{aodd}
a_{1_{1},\cdots , i_{2k+1}}= \mathbf{E}\left[ \frac{1}{V_{n}}W(\varphi_{1} ) \hp_{2d_{1}+1} \left(  W(\varphi_{1})\right) \hp_{2d_{2}} \left(  W(\varphi_{2})\right)\cdots  \hp_{2d_{n}} \left(  W(\varphi_{n})\right)\right].
\end{equation}
\end{remark}
\noindent {\bf Proof: }Since $\sum_{r=1}^{n}d_r^{\star} = 2k+1$, there is an odd number of odd integers in the sequence $d_{r}^{\star}, 1 \leq r \leq n$. Recall that by Lemma \ref{lem:intbypart}, we have
\begin{eqnarray}
\label{ecritDev}
a_{i_{1}, \cdots  i_{2k+1}} &=& \sum_{u=1}^{n} \mathbf{E} \left( \frac{W(\varphi_u)}{V_{n}} \prod_{r=1}^n \hp_{d_r^{\star}}\big(W(\varphi _{r})\big)  \right)\nonumber \\
&=& \mathbf{E}\left[ \frac{1}{V_{n}}W(\varphi_{1} ) \hp_{d_{1}^{\star}} \left(  W(\varphi_{1})\right) \hp_{d_{2}^{\star}} \left(  W(\varphi_{2})\right)\ldots \hp_{d_{n}^{\star}} \left(  W(\varphi_{n})\right)\right]\nonumber \\
&+& \mathbf{E}\left[ \frac{1}{V_{n}}W(\varphi_{2} ) \hp_{d_{1}^{\star}} \left(  W(\varphi_{1})\right) \hp_{d_{2}^{\star}} \left(  W(\varphi_{2})\right)\cdots  \hp_{d_{n}^{\star}} \left(  W(\varphi_{n})\right)\right]\nonumber \\
&\vdots &  \nonumber \\
&+& \mathbf{E}\left[ \frac{1}{V_{n}}W(\varphi_{n} ) \hp_{d_{1}^{\star}} \left(  W(\varphi_{1})\right) \hp_{d_{2}^{\star}} \left(  W(\varphi_{2})\right)\cdots  \hp_{d_{n}^{\star}} \left(  W(\varphi_{n})\right)\right].
\end{eqnarray}
Because of Lemma \ref{lemmaxx}, for each $i$, the term $$\mathbf{E}\left[ \frac{1}{V_{n}}W(\varphi_{i} ) \hp_{d_{1}^{\star}} \left(  W(\varphi_{1})\right) \hp_{d_{2}^{\star}} \left(  W(\varphi_{2})\right)\cdots  \hp_{d_{n}^{\star}} \left(  W(\varphi_{n})\right)\right]$$ is non null if and only if $d_{i}^{\star}$ is the only odd integer in $d_{r}^{\star}, 1 \leq r \leq n$. Thus, $a_{1_{1},\cdots , i_{2k+1}} \neq 0$ if there is only one odd integer in $d_{r}^{\star}, 1 \leq r \leq n$. Let $d_{i}^{\star}$ with $1\leq i \leq n$ be this only odd integer. Then, if $j \neq i$, by Lemma \ref{lemmaxx},
$$\mathbf{E}\left[ \frac{1}{V_{n}}W(\varphi_{j} ) \hp_{d_{1}^{\star}} \left(  W(\varphi_{1})\right) \hp_{d_{2}^{\star}} \left(  W(\varphi_{2})\right)\cdots  \hp_{d_{n}^{\star}} \left(  W(\varphi_{n})\right)\right] = 0.$$ Thus, using (\ref{ecritDev}) yields
\begin{equation*}
a_{1_{1},\cdots , i_{2k+1}}= \mathbf{E}\left[ \frac{1}{V_{n}}W(\varphi_{i} ) \hp_{d_{1}^{\star}} \left(  W(\varphi_{1})\right) \hp_{d_{2}^{\star}} \left(  W(\varphi_{2})\right)\cdots  \hp_{d_{n}^{\star}} \left(  W(\varphi_{n})\right)\right]
\end{equation*}
if there is only one odd integer in the sequence $d_{r}^{\star}, 1 \leq r \leq n$ and $a_{1_{1},\cdots , i_{2k+1}}= 0$ if there is more than one odd integer in the sequence $d_{r}^{\star}, 1 \leq r \leq n$. Using the equality in law between $W(\varphi_{i} )$ and $W(\varphi_{1} )$, one can perform an index swap $(i\leftrightarrow 1)$ to finally obtain the desired result.
\qed \vskip0.3cm

\noindent In the following lemma, we compute the $L^{2}$ norm of $F_n$. This technical result will be needed in the next section.
\begin{lemma} \label{lem:Chaosnorm} Let $a_{i_{1},\cdots , i_{2k+1}}$ be as given in \eqref{eqn:chaosFN1}. Then, for every $n \in \mathbb{N}$, we have
\begin{equation*}
\left\| F_n \right\|_{L^{2}(\Omega)}^{2} = \sum_{k\geq 0} \frac{1}{(2k+1)! } \sum_{i_{1},\cdots , i_{2k+1}  = 1}^n a_{i_{1},\cdots ,i_{2k+1}}^{2} = 1.
\end{equation*}
\end{lemma}
\noindent {\bf Proof: }Firstly, using the isometry of multiple stochastic integrals and the orthogonality of the kernels $\varphi_{i}$, one can write
\begin{eqnarray*}
\mathbf{E} \left( F_{n}^{2}\right)  &=& \sum_{k\geq 0} \left( \frac{1}{(2k+1)! }\right) ^{2}(2k+1)!\sum_{ \substack{i_{1},\cdots , i_{2k+1}=1 \\ j_{1},\cdots , j_{2k+1}=1}}
^{n}  a_{i_{1},\cdots , i_{2k+1}} a_{j_{1},\cdots , j_{2k+1}} \nonumber \\
&\times &  \Big \langle   \varphi_{i_{1}}\otimes  \cdots  \otimes \varphi_{i_{2k+1}},  \varphi_{j_{1}}\otimes \cdots  \otimes \varphi_{j_{2k+1}} \Big \rangle _{L^{2}([0,1] ^{2k})}\nonumber  \\
&=&\sum_{k\geq 0} \frac{1}{(2k+1)! } \sum_{i_{1},\cdots , i_{2k+1}=1}^n a_{i_{1},\cdots ,i_{2k+1}}^{2}.
\end{eqnarray*}
Secondly, using the fact that $F_n ^{2} =\frac{W_{n}^{2}}{V_{n}^{2}}$, we have
\begin{eqnarray*}
\mathbf{E}\left( F_{n}^{2}\right) &=& \frac{1}{(2\pi)^{\frac{n}{2}}}\int_{\mathbb{R}^{n}}\frac{(x_{1}+\cdots +x_{n})^{2} }{x_{1}^{2}+\cdots  +x_{n}^{2}}e^{-\frac{1}{2}(x_{1}^{2}+\cdots  +x_{n}^{2})}dx_{1}\cdots  dx_{n}\\
&=&\frac{1}{(2\pi)^{\frac{n}{2}}}\int_{\mathbb{R}^{n}}\frac{x_{1}^{2}+\cdots  +x_{n}^{2} }{x_{1}^{2}+\cdots  +x_{n}^{2}}e^{-\frac{1}{2}(x_{1}^{2}+\cdots  +x_{n}^{2})}dx_{1}\cdots  dx_{n}=1
\end{eqnarray*}
because the mixed terms vanish as in the proof of Lemma \ref{lem:evenzero}. \qed \vskip0.3cm

\noindent Recall that if $X$ is a Chi-squared random variable with $n$ degrees of
freedom (denoted by $\chi _{n} ^{2}$)  then for any $m \geq 0$,
\begin{equation*}
\mathbf{E} \left( X^{m}\right) = 2^{m} \frac{ \Gamma (m+\frac{n}{2})}{\Gamma (\frac{n}{2})}.
\end{equation*}
where $\Ga(\cdot)$ denotes the standard Gamma function.
\\~\\
\noindent When $k=0$, the coefficients $a_{i_{1},\cdots ,i_{2k+1}}$ can be easily computed. Indeed,
noticing that $V_n^2$ has a  $\chi_n^2$ distribution, we obtain
\begin{equation*}
\sum_{i=1}^{n} a_{i} = \mathbf{E} \left( \sum_{i=1}^{n} \frac{1}{V_{n}}W(\varphi_{i})^{2}\right)  = \mathbf{E} \left( (V_{n}^{2} ) ^{\frac{1}{2}} \right) =2^{\frac{1}{2}}\frac{ \Gamma (\frac{1}{2}+\frac{n}{2})}{\Gamma (\frac{n}{2})}.
\end{equation*}

\noindent Since $a_{1} =a_{2}=\cdots =a_{n}$ we obtain that for every $i=1,..,n$
\begin{equation*}
a_{i} = \frac{2^{\frac{1}{2}}}{n}\frac{ \Gamma (\frac{1}{2}+\frac{n}{2})}{\Gamma (\frac{n}{2})}.
\end{equation*}

\vskip0.2cm

\noindent The following lemma is the second key result in our goal of calculating the coefficients. It will be used repeatedly in the sequel.
\begin{lemma}\label{lem:main}  Let $\{a_1, a_2,\cdots a_n\} $ be
non-negative numbers. Then it holds that
\begin{eqnarray*}
&\mathbf{E} \left ({W(\vp_1)^{2a_1} W(\vp_2)^{2a_2} \cdots W(\vp_n)^{2a_n} \over V_n }\right )\nonumber\\& = {1 \over (2\pi)^{\frac{n}{2}}}2^{a_{1}+\cdots +a_{n}+\frac{n-1}{2}}\frac{ \Gamma \left(a_{1}+\cdots +a_{n}+\frac{n-1}{2} \right)}{\Gamma \left(a_{1}+\cdots +a_{n}+\frac{n}{2} \right)}\Gamma\left(a_{1}+\frac{1}{2}\right)\cdots  \Gamma\left(a_{n}+\frac{1}{2}\right).
\end{eqnarray*}
\end{lemma}
\noindent {\bf Proof: }By definition, we have
\begin{eqnarray*}
&& \mathbf{E} \left ({W(\vp_1)^{2a_1} W(\vp_2)^{2a_2} \cdots W(\vp_n)^{2a_n} \over V_n }\right )\nonumber \\
&& = \frac{1}{(2\pi)^{\frac{n}{2}}}\int_{\mathbb{R}^n}{x_1^{2a_1} x_2^{2a_2}\cdots x_n^{2a_n} \over \sqrt{x_1^2 +
x_2^2 + \cdots + x_n^2}} e^{-{1\over 2}(x_1^2 + x_2^2 + \cdots + x_n^2)} dx_1 dx_2 \cdots  dx_n\nonumber \\
&& = {1 \over (2\pi)^{\frac{n}{2}}}I.
\end{eqnarray*}
To compute the above integral $I$, we introduce $n$-dimensional polar coordinates.
Set
\begin{eqnarray*}
&& x_1 = r \cos \theta_1 \\
&& x_j = r \cos \theta_j \prod_{r=1}^{j-1} \sin \theta_r,  \quad \,  j= 2, \cdots , n-2 \\
&& x_{n-1} = r \sin \psi \prod_{r=1}^{n-2} \sin \theta_r , \quad \, x_n = r \cos \psi \prod_{r=1}^{n-2} \sin \theta_r
\end{eqnarray*}
with $ 0 \leq r < \infty $, $0 \leq \theta_i \leq \pi$ and $0 \leq \psi \leq 2\pi$. It can be easily verified
that $x_1^2 + x_2^2 + \cdots + x_n^2 = r^2$.
The Jacobian of the above transformation is given by
$$J = r^{n-1} \prod_{k=1}^{n-2} \sin^{k} \theta_{n-1-k} \,.$$
Therefore our integral denoted by $I$ becomes
\begin{eqnarray*}
&&\int_{0}^{\infty} r^{2(a_{1}+\cdots + a_{n}) +n-2} e^{-\frac{r^{2}}{2}}dr \int_{0}^{2\pi} (\sin \psi ) ^{2a_{n-1} + 2a_{n}} (\cos \psi ) ^{2a_{n} }d\psi\\
&&\prod _{k=2}^{n-1}\int_{0}^{\pi}(\sin \theta _{n-k} ) ^{2a_{n}+2a_{n-1}+\cdots + 2a_{n-k+1}+k-1} (\cos \theta _{n-k}) ^{2a_{n-k}}d\theta _{n-k}.
\end{eqnarray*}
Let us compute the first integral with respect to $dr$. Using the change of variables $\frac{r^{2}}{2}=y$, we get
\begin{eqnarray*}
\int_{0}^{\infty} r^{2(a_{1}+\cdots + a_{n}) +n-2} e^{-\frac{r^{2}}{2}}dr&=&
2^{a_{1}+\cdots +a_{n}+\frac{n-1}{2}-1}\int_{0}^{\infty} dy y^{a_{1}+\cdots +a_{n} +\frac{n-1}{2}-1 }e^{-y}\\
&=& 2^{a_{1}+\cdots +a_{n}+\frac{n-1}{2}-1}\Gamma \left( a_{1}+\cdots +a_{n}+\frac{n-1}{2}\right).
\end{eqnarray*}
Let us now compute the integral with respect to $d\psi$. We use the following formula: for every $a,b\in \mathbb{Z}$, it holds that
\begin{eqnarray*}
\int_{0}^{2\pi} (\sin \theta ) ^{a} (\cos \theta )^{b}d\theta &=
2 \beta \left(\frac{a+1}{2}, \frac{b+1}{2} \right) \mbox{if $m$ and $n$ are even}\\
&= 0, \hskip0.5cm \mbox{ if $m$ or $n$ are odd. }
\end{eqnarray*}
This implies that
\begin{equation*}
\int_{0}^{2\pi} (\sin \psi ) ^{2a_{n-1} + 2a_{n}} (\cos \psi ) ^{2a_{n} }d\psi=2\beta \left( a_{n}+\frac{1}{2}, a_{n-1}+\frac{1}{2}\right).
\end{equation*}
Finally, we deal with the integral with respect to $d\theta _{i}$  for $i=1$ to $n-2$. Using the fact that, for $a,b>-1$, it holds that
\begin{equation*}
\int_{0}^{\frac{\pi}{2}} (\sin \theta ) ^{a} (\cos \theta )^{b}d\theta= \frac{1}{2} \beta \left(\frac{a+1}{2}, \frac{b+1}{2} \right)
\end{equation*}
yields
\begin{eqnarray*}
&&\int_{0}^{\pi} (\sin \theta _{n-k} ) ^{2a_{n}+2a_{n-1}+\cdots + 2a_{n-k+1}+k-1} (\cos \theta _{n-k}) ^{2a_{n-k}}d\theta _{n-k}\\
&=&
\int_{0}^{\frac{\pi}{2}}(\sin \theta _{n-k} ) ^{2a_{n}+2a_{n-1}+\cdots + 2a_{n-k+1}+k-1} (\cos \theta _{n-k}) ^{2a_{n-k}}d\theta _{n-k}\\
&&+
\int_{\frac{\pi}{2}}^{\pi} (\sin \theta _{n-k} ) ^{2a_{n}+2a_{n-1}+\cdots + 2a_{n-k+1}+k-1} (\cos \theta _{n-k}) ^{2a_{n-k}}d\theta _{n-k}\\
&=& \frac{1}{2} \beta \left( a_{n}+\cdots +a_{n-k+1}+ \frac{k}{2}, a_{n-k}+\frac{1}{2}\right)\\
&&+ \int_{0}^{\frac{\pi}{2}}(\sin (\theta _{n-k}+\frac{\pi}{2}) ) ^{2a_{n}+2a_{n-1}+\cdots + 2a_{n-k+1}+k-1} (\cos (\theta _{n-k}+\frac{\pi}{2})) ^{2a_{n-k}}d\theta _{n-k}\\
&=&\beta \left( a_{n}+\cdots +a_{n-k+1}+ \frac{k}{2}, a_{n-k}+\frac{1}{2}\right)
\end{eqnarray*}
because $\sin(\theta +\frac{\pi}{2})= \cos \theta $ and $\cos (\theta + \frac{\pi}{2})=-\sin (\theta)$. By gathering the above calculations, the integral $I$ becomes
\begin{eqnarray*}
I&=&2^{a_{1}+\cdots +a_{n}+\frac{n-1}{2}}\Gamma \left(a_{1}+\cdots +a_{n}+\frac{n-1}{2} \right) \beta \left( a_{n}+\frac{1}{2}, a_{n-1}+\frac{1}{2}\right)\\
&&\times \prod_{k=2}^{n-1} \beta \left( a_{n}+\cdots +a_{n-k+1}+ \frac{k}{2}, a_{n-k}+\frac{1}{2}\right)\\
&=&2^{a_{1}+\cdots +a_{n}+\frac{n-1}{2}}\Gamma \left(a_{1}+\cdots +a_{n}+\frac{n-1}{2} \right) \frac{\Gamma \left( a_{n}+\frac{1}{2}\right)\Gamma \left( a_{n-1}+\frac{1}{2}\right)}{\Gamma \left( a_{n}+a_{n-1}+1\right) } \\
&&\times \prod_{k=2}^{n-1}\frac{ \Gamma  \left( a_{n}+\cdots +a_{n-k+1}+ \frac{k}{2}\right) \Gamma \left(a_{n-k}+\frac{1}{2}\right)}{\Gamma \left( a_{n}+a_{n-1}+\cdots + a_{n-k}+\frac{k+1}{2}\right)}\\
&=&2^{a_{1}+\cdots +a_{n}+\frac{n-1}{2}}\frac{ \Gamma \left(a_{1}+\cdots +a_{n}+\frac{n-1}{2} \right)}{\Gamma \left(a_{1}+\cdots +a_{n}+\frac{n}{2} \right)}\Gamma\left(a_{1}+\frac{1}{2}\right)\cdots  \Gamma\left(a_{n}+\frac{1}{2}\right).
\end{eqnarray*}
This concludes the proof. \qed \vskip0.3cm

\noindent The following result gives the asymptotic behavior of the coefficients when $n\to \infty$.
\begin{lemma}\label{borne-coef}
For every $1\leq i_{1},\cdots , i_{2k+1} \leq n$, let $a_{i_{1}, \cdots , i_{2k+1}}$
be as defined in \eqref{eqn:ai}. As in (\ref{aodd}), let $2d_{1}+1, 2d_{2},\cdots , 2d_r,\cdots ,2d_{n}$ denote
the number of times the integer $r$ appears in the sequence $\{i_1,i_2,\cdots , i_{2k+1}\}$ with
$\sum_{r=1}^n d_r = k$. Then when $n\to \infty$,
\begin{eqnarray}
\label{imp}
a_{i_{1},\cdots , i_{2k+1}} & \sim & \frac{1}{k!} (2k-1)!!\frac{(2d_{1}+1)! (2d_{2})! \cdots (2d_{n})! }{(d_{1}! d_{2}! \cdots d_{n}!)^{2}} \nonumber \\
&&  \times \ 2^{-2k}(-1)^{k}\left( \prod _{j=0}^{n} \sum_{l_{j}=0}^{d_{j}}(-1)^{l_{j}}C_{d_{j}}^{l_{j}}l_{j} ^{d_{j}}\right)\frac{1}{n^{\frac{1}{2}+\left| A\right|}}
\end{eqnarray}
where
\begin{equation*}
A:=  \{ 2d_{1} +1, 2d_{2},\cdots ,2d_{n}\} \setminus \{0,1 \}
\end{equation*}
and $\left| A\right|$ is the cardinal of $A$.
\end{lemma}
\noindent {\bf Proof: }We recall the following explicit formula for the Hermite polynomials
\begin{equation}\label{exp}
\hp_{d}(x)=d! \sum_{l=0}^{[\frac{d}{2}] } \frac{ (-1)^{l} }{ 2^{l} l! (d-2l)! } x^{d-2l}.
\end{equation}

\noindent Using (\ref{exp}) and (\ref{aodd}) we can write
\begin{eqnarray*}
&&a_{i_{1},\cdots , i_{2k+1}}=\mathbf{E}\left[ \frac{1}{V_{n}}W(\varphi_{1} ) \hp_{2d_{1} +1} \left(  W(\varphi_{1})\right) \hp_{2d_{2}} \left(  W(\varphi_{2})\right)\cdots \hp_{2d_{n} } \left(  W(\varphi_{n})\right)\right]\\
&=&(2d_{1}+1)! (2d_{2})! \cdots (2d_{n})! \sum_{l_{1}=0}^{d_{1}}\sum_{l_{2}=0}^{d_{2}}\cdots \sum_{l_{n}=0}^{d_{n}}\frac{(-1)^{l_{1}+l_{2}+\cdots  +l_{n}}}{2^{l_{1}+l_{2}+\cdots + l_{n}}l_{1}!\cdots l_{n}!}\\
&&\times\ \frac{\mathbf{E} \left[\frac{1}{V_{n}} W(\varphi_{1} ) ^{2d_{1}+2-2l_{1}}  W(\varphi_{2} ) ^{2d_{2}-2l_{2}}\cdots W(\varphi_{n} ) ^{2d_{n}-2l_{n}}\right]}{(2d_{1}+1-2l_{2})! (2d_{2}-2l_{2})! \cdots (2d_{n}-2l_{n})!}.
\end{eqnarray*}
At this point, we use Lemma \ref{lem:main} to rewrite the expectation in the last equation.

\begin{eqnarray*}
&& \mathbf{E}\left[ \frac{1}{V_{n}}W(\varphi_{1} ) \hp_{2d_{1} +1} \left(  W(\varphi_{1})\right) \hp_{2d_{2}} \left(  W(\varphi_{2})\right)\cdots  \hp_{2d_{n} } \left(  W(\varphi_{n})\right)\right]\\
&=&(2d_{1}+1)! (2d_{2})! \cdots (2d_{n})! \sum_{l_{1}=0}^{d_{1}}\sum_{l_{2}=0}^{d_{2}}\cdots  \sum_{l_{n}=0}^{d_{n}}\frac{(-1)^{l_{1}+l_{2}+\cdots +l_{n}}}{2^{l_{1}+l_{2}+\cdots + l_{n}}l_{1}!\cdots l_{n}!}\\
&& \times \frac{2^{d_{1}+1+d_{2}+\cdots +d_{n}-(l_{1}+l_{2}+\cdots +l_{n})+\frac{n-1}{2}}}{(2\pi)^{\frac{n}{2}}(2d_{1}+1-2l_{2})! (2d_{2}-2l_{2})! \cdots (2d_{n}-2l_{n})!} \\
&& \times \frac{ \Gamma \left(d_{1}+1+d_{2}+\cdots +d_{n}-(l_{1}+l_{2}+\cdots +l_{n})+\frac{n-1}{2} \right)}{\Gamma \left(d_{1}+1+d_{2}+\cdots +d_{n}-(l_{1}+l_{2}+\cdots +l_{n})+\frac{n}{2} \right)}\\
&& \times \Gamma\left(d_{1}+1-l_{1}+\frac{1}{2}\right)\Gamma\left(d_{2}-l_{2}+\frac{1}{2}\right)\cdots  \Gamma\left(d_{n}-l_{n}+\frac{1}{2}\right)\\
&=&(2d_{1}+1)! (2d_{2})! \cdots (2d_{n})! \sum_{l_{1}=0}^{d_{1}}\sum_{l_{2}=0}^{d_{2}}\cdots  \sum_{l_{n}=0}^{d_{n}}\frac{(-1)^{l_{1}+l_{2}+\cdots +l_{n}}}{2^{2(l_{1}+l_{2}+\cdots + l_{n})}l_{1}!\cdots l_{n}! }\\
&&\times \frac{2^{d_{1}+1+d_{2}+\cdots +d_{n}-\frac{1}{2}}}{\pi^{\frac{n}{2}}(2d_{1}+1-2l_{2})! (2d_{2}-2l_{2})! \cdots (2d_{n}-2l_{n})!}\\
&& \times \frac{ \Gamma \left(d_{1}+1+d_{2}+\cdots +d_{n}-(l_{1}+l_{2}+\cdots +l_{n})+\frac{n-1}{2} \right)}{\Gamma \left(d_{1}+1+d_{2}+\cdots +d_{n}-(l_{1}+l_{2}+\cdots +l_{n})+\frac{n}{2} \right)}\\
&&\times \Gamma\left(d_{1}+1-l_{1}+\frac{1}{2}\right)\Gamma\left(d_{2}-l_{2}+\frac{1}{2}\right)\cdots  \Gamma\left(d_{n}-l_{n}+\frac{1}{2}\right).
\end{eqnarray*}
\noindent We claim that for any integers $d\geq l$,
\begin{equation}
\label{claim1}
\frac{(-1)^{l}}{2^{-2l}l! (2d-2l)! }\Gamma \left(d-l+\frac{1}{2} \right)=\sqrt{\pi} \frac{ 2^{-2d}(-1)^{l}}{d!}C_{d}^{l}.
\end{equation}
Recall the relation satisfied by the Gamma function: for every $z>0$,
\begin{equation}
\label{relationgamma}
 \Gamma (z+1)=z\Gamma (z) \mbox{ and } \Gamma (z) \Gamma (z+\frac{1}{2})= \sqrt{\pi} 2^{1-2z}\Gamma (2z).
 \end{equation}
Then
\begin{eqnarray*}
\frac{(-1)^{l}}{2^{-2l}l! (2d-2l)! }\Gamma \left(d-l+\frac{1}{2} \right)&=&\frac{(-1)^{l}}{2^{-2l}l! (2d-2l)! }\frac{\Gamma\left( d-l+1+\frac{1}{2}\right)}{d-l-\frac{1}{2}}\\
&=& \frac{(-1)^{l}}{2^{-2l}l! (2d-2l)! }\frac{\Gamma (2d-2l+2) }{\Gamma(d-l+1)} \sqrt{\pi} 2^{1-2(d-l+1)}\\
&=& \sqrt{\pi} 2^{-2d}\frac{(-1)^{l}}{l! (2d-2l)! } \frac{(2d-2l+1)! }{(d-l)! (2d-2l+1)}\\&=&\sqrt{\pi} \frac{ 2^{-2d}(-1)^{l}}{d!}C_{d}^{l}
\end{eqnarray*}
and (\ref{claim1}) is proved. In the same way, using only the second relation in (\ref{relationgamma}), we obtain
\begin{equation}\label{claim2}
\frac{(-1)^{l_{1}}}{2^{-2l_{1}}l_{1}! (2d_{1}+1-2l_{1})! }\Gamma \left(d_{1}+1-l_{1}+\frac{1}{2} \right)=\sqrt{\pi} \frac{ 2^{-1-2d_{1}}(-1)^{l_{1}}}{d_{1}!}C_{d_{1}}^{l_{1}}.
\end{equation}
Putting together (\ref{claim1}) and (\ref{claim2}) we find
\begin{eqnarray*}
&& \mathbf{E}\left[ \frac{1}{V_{n}}W(\varphi_{1} ) \hp_{2d_{1} +1} \left(  W(\varphi_{1})\right) \hp_{2d_{2}} \left(  W(\varphi_{2})\right)\cdots  \hp_{2d_{n} } \left(  W(\varphi_{n})\right)\right]\\
&=&\frac{(2d_{1}+1)! (2d_{2})! \cdots (2d_{n})! }{d_{1}! d_{2}! \cdots d_{n}!}2^{-(d_{1}+\cdots +d_{n})-\frac{1}{2}}\sum_{l_{1}=0}^{d_{1}}\sum_{l_{2}=0}^{d_{2}}\cdots  \sum_{l_{n}=0}^{d_{n}}
(-1)^{l_{1}+l_{2}+\cdots +l_{n}}C_{d_{1}}^{l_{1}}\cdots C_{d_{n}}^{l_{n}}\\
&&\times \frac{ \Gamma \left(d_{1}+1+d_{2}+\cdots +d_{n}-(l_{1}+l_{2}+\cdots +l_{n})+\frac{n-1}{2} \right)}{\Gamma \left(d_{1}+1+d_{2}+\cdots +d_{n}-(l_{1}+l_{2}+\cdots +l_{n})+\frac{n}{2} \right)}.
\end{eqnarray*}
By Stirling's formula, when $n$ goes to infinity, we have
\begin{equation*}
\frac{ \Gamma \left(d_{1}+1+d_{2}+\cdots +d_{n}-(l_{1}+l_{2}+\cdots +l_{n})+\frac{n-1}{2} \right)}{\Gamma \left(d_{1}+1+d_{2}+\cdots +d_{n}-(l_{1}+l_{2}+\cdots +l_{n})+\frac{n}{2} \right)}\sim \frac{1}{\sqrt{ k+1-(l_{1}+\cdots +l_{n})+ \frac{n}{2}}}.
\end{equation*}
Therefore we need to study the behavior of the sequence
$$t_{n}:=\sum_{l_{1}=0}^{d_{1}}\sum_{l_{2}=0}^{d_{2}}\cdots  \sum_{l_{n}=0}^{d_{n}}
(-1)^{l_{1}+l_{2}+\cdots +l_{n}}C_{d_{1}}^{l_{1}}\cdots C_{d_{n}}^{l_{n}}\frac{1}{\sqrt{ k+1-(l_{1}+\cdots +l_{n})+ \frac{n}{2}}}$$
as $n\to \infty$.
We can write
$$t_{n}=\frac{1}{\sqrt{n}}\sqrt{2}g(\frac{1}{n}) $$
where
\begin{equation*}
g(x) = \sum_{l_{1}=0}^{d_{1}}\sum_{l_{2}=0}^{d_{2}}\cdots  \sum_{l_{n}=0}^{d_{n}}
(-1)^{l_{1}+l_{2}+\cdots +l_{n}}C_{d_{1}}^{l_{1}}\cdots C_{d_{n}}^{l_{n}}\frac{1}{\sqrt{2k+2-(l_{1}+\cdots +l_{n}) x +1} } .
\end{equation*}
Since for every $d\geq 1$
\begin{equation*}
\sum_{l=0}^{d}(-1)^{l}C_{d}^{l}=0
\end{equation*}
we clearly have $g(0)=0$. The $q^{\mbox{\tiny{th}}}$ derivative of $g$ at zero is
\begin{equation*}
g^{(q)}(0)= (-1)^{q} \frac{(2q-1)!!}{2^{q}} \left[ 2k+2-(l_{1}+\cdots +l_{n})\right] ^{q}.
\end{equation*}
Repeatedly using the relation $C_{n}^{k} =\frac{n}{k}C_{n-1}^{k-1}$ we can prove  that
\begin{equation*}
\sum_{l=0}^{d}(-1)^{l}C_{d}^{l}l^{q}=0
\end{equation*}for every $q=0,1,\cdots , d-1$.
Therefore the first non-zero term in the Taylor decomposition of the function $g$ around zero is $$  \sum_{l_{1}=0}^{d_{1}}\sum_{l_{2}=0}^{d_{2}}\cdots  \sum_{l_{n}=0}^{d_{n}}
(-1)^{l_{1}+l_{2}+\cdots +l_{n}}C_{d_{1}}^{l_{1}}\cdots C_{d_{n}}^{l_{n}}l_{1}^{d_{1}}\cdots l_{n}^{d_{n}}$$ which appears when we take the derivative of order $d_{1}+d_{2}+\cdots +d_{n}$.
We obtain that, for $x$ close to zero,
\begin{equation*}
g(x)\sim (-1)^{d_{1}+\cdots +d_{n}} \frac{(2(d_{1}+\cdots +d_{n})-1)!!}{2^{d_{1}+\cdots +d_{n}}}\prod _{j=0}^{n} \sum_{l_{j}=0}^{d_{j}}(-1)^{l_{j}}C_{d_{j}}^{l_{j}}l_{j} ^{d_{j}}\times H(d_{1},\cdots , d_{n})x^{\left| A\right| }
\end{equation*}
where
$$A=\{ d_{1},\cdots , d_{n}\} \setminus \{0\}= \{ 2d_{1} +1, 2d_{2},\cdots ,2d_{n}\}\setminus \{0,1\} $$
and $H(d_{1},\cdots , d_{n})$ is the coefficient of $l_{1}^{d_{1}} \cdots l_{n}^{d_{n}}$ in the expansion of $(l_{1}+\cdots + l_{n}) ^{d_{1}+\cdots +d_{n}}$. That is
\begin{equation*}
H(d_{1},\cdots , d_{n})=C_{d_{1}+\cdots +d_{n}}^{d_{1}}C_{d_{2}+\cdots +d_{n}}^{d_{2}}\cdots C_{d_{n-1}+d_{n}}^{d_{n-1}}= \frac{(d_{1}+\cdots +d_{n})! }{d_{1}!\cdots d_{n}!} .
\end{equation*}
We finally have
\begin{eqnarray*}
a_{i_{1},\cdots ,i_{2k+1}} &=& \frac{(2d_{1}+1)! (2d_{2})! \cdots (2d_{n})! }{(d_{1}! d_{2}! \cdots d_{n}!)^{2}}2^{-(d_{1}+\cdots +d_{n})}(-1)^{d_{1}+\cdots +d_{n}} \frac{(2(d_{1}+\cdots +d_{n})-1)!!}{2^{d_{1}+\cdots +d_{n}}}\\
&&\times
\left( \prod _{j=0}^{n} \sum_{l_{j}=0}^{d_{j}}(-1)^{l_{j}}C_{d_{j}}^{l_{j}}l_{j} ^{d_{j}}\right)\frac{(d_{1}+\cdots +d_{n})! }{d_{1}!\cdots d_{n}!}\frac{1}{n^{\frac{1}{2}+\left| A\right|}}\\
&=&k! (2k-1)!!\frac{(2d_{1}+1)! (2d_{2})! \cdots (2d_{n})! }{(d_{1}! d_{2}! \cdots d_{n}!)^{2}}2^{-2k}(-1)^{k}\left( \prod _{j=0}^{n} \sum_{l_{j}=0}^{d_{j}}(-1)^{l_{j}}C_{d_{j}}^{l_{j}}l_{j} ^{d_{j}}\right)\frac{1}{n^{\frac{1}{2}+\left| A\right|}}\\
&=& k! (2k-1)!!\frac{(2d_{1}+1)! (2d_{2})! \cdots (2d_{n})! }{(d_{1}! d_{2}! \cdots d_{n}!)^{2}}2^{-2k}(-1)^{k}\left( \prod _{j=0}^{n} t(d_{j})\right)\frac{1}{n^{\frac{1}{2}+\left| A\right|}}
\end{eqnarray*}
with for $i=1,\cdots ,n$
\begin{equation}
\label{tdj}
t(d_{j}):=\sum_{l_{j}=0}^{d_{j}}(-1)^{l_{j}}C_{d_{j}}^{l_{j}}l_{j} ^{d_{j}}.
\end{equation}
\qed \vskip0.3cm

\section{Computation of the Berry-Ess\'een bound}

\noindent Let us first recall the following result (see \cite{DLS}, page 53): for $2<p\leq 3$,
\begin{equation}
\label{borneDLSbis}
\sup_{z\in \mathbb{R}} \left| P(F_{n} \leq z) -\Phi (z) \right| \leq 25 \mathbf{E}\left( |Z| ^{p}\right)  n^{1-\frac{p}{2} }
\end{equation}
where $Z$ is a standard normal random variable and $\Phi$ is its repartition function. In particular for $p=3$ we get
$$\sup_{z\in \mathbb{R}} \left| P(F_{n} \leq z) -\Phi (z) \right| \leq 25 \mathbf{E}\left( |Z| ^{3}\right) n^{-\frac{1}{2} }.$$

\noindent We now compute the Berry-Essen bound obtained via Malliavin calculus in order to compare it  with (\ref{borneDLSbis}). Formula (\ref{eqn:chaosFN1}) yields
\begin{equation}
\label{deriveMdeFn}
D_{\alpha }F_{n}= \sum_{k\geq 0} \frac{2k+1}{(2k+1)!}  \sum_{i_{1}, \cdots , i_{2k+1}=1}^{n} a_{i_{1},\cdots , i_{2k+1}} I_{2k} \left( (\varphi_{i_{1}}\otimes  \cdots \otimes
\varphi_{i_{2k=1}})^{\sim} \right)(\cdot , \alpha)
\end{equation}
(here $(\varphi_{i_{1}}\otimes  \cdots \otimes \varphi_{i_{2k=1}})^{\sim}$ denotes the symmetrization of the function $\varphi_{i_{1}}\otimes  \cdots \otimes \varphi_{i_{k}}$ with respect to its $k$ variables) and
\begin{equation}
\label{deriveMLmoins1deFn}
D_{\alpha} (-L)^{-1}F_{n} = \sum_{k}\frac{1}{(2k+1)!}\sum_{i_{1}, \cdots , i_{2k+1}=1}^{n} a_{i_{1},\cdots , i_{2k+1}} I_{2k} \left( ( \varphi_{i_{1}}\otimes  \cdots \otimes \varphi_{i_{2k+1}} )^{\sim }\right) (\cdot , \alpha ).
\end{equation}

\noindent It is now possible to calculate the quantity
\begin{equation*}
\mathbf{E}\left(  1- \langle DF_{n}, D(-L)^{-1}F_{n} \rangle \right) ^{2}
\end{equation*}
more explicitly by using the product formula (\ref{prod}) and the isometry of multiple stochastic integrals.
\begin{lemma}
For every $n\geq 2$,
\begin{eqnarray*}
&&\mathbf{E}\left( 1- \langle DF_{n}, D(-L)^{-1}F_{n} \rangle \right) ^{2}\\
&=&\sum_{m\geq 1}(2m)! \sum_{i_{1},\cdots , i_{2m} =1}^{n}\left( \sum_{k=0}^{2m} \frac{1}{k!} \frac{1}{(2m-k)!}\sum_{r\geq 0} \frac{1}{r!}\frac{1}{2m-k+r+1}
\sum_{u_{1}, \cdots , u_{r+1}=1}^{n}\right. \\
&&\left. a_{u_{1}, u_{2}, \cdots , u_{r+1}, i_{1}, \cdots , i_{k}} a_{u_{1}, u_{2}, \cdots , u_{r+1}, i_{k+1}, \cdots , i_{2m}}\Bigg )^{2}\right. .
\end{eqnarray*}
\end{lemma}
\noindent {\bf Proof: }Using (\ref{deriveMdeFn}) and (\ref{deriveMLmoins1deFn}), we can calculate the following quantity.
\begin{eqnarray*}
 \langle DF_{n}, D(-L)^{-1}F_{n} \rangle &=& \sum_{k,l\geq 0} \frac{1}{(2k)! }\frac{1}{(2l+1)! } \sum_{i_{1}, \cdots , i_{2k+1}=1}^{n} a_{i_{1},\cdots , i_{2k+1}}\sum_{j_{1}, \cdots , j_{2l+1}=1}^{n} a_{j_{1},\cdots , j_{2l+1}}\\
&&\times \int_{0}^{\infty} d\alpha I_{2k}  \left( (\varphi_{i_{1}}\otimes  \cdots \otimes
\varphi_{i_{2k+1}})^{\sim} \right)(\cdot , \alpha)I_{2l}  \left( (\varphi_{j_{1}}\otimes  \cdots \otimes
\varphi_{j_{2l+1}})^{\sim} \right)(\cdot , \alpha)\\
&=& \sum_{k,l\geq 0} \frac{1}{(2k)! }\frac{1}{(2l+1)! }\sum_{u=1}^{n} \sum_{i_{1}, \cdots , i_{2k}=1}^{n} a_{u, i_{1},\cdots , i_{2k}}\sum_{j_{1}, \cdots , j_{2l}=1}^{n} a_{u,j_{1},\cdots ,
j_{2l}}\\
&&\times  I_{2k}  \left( (\varphi_{i_{1}}\otimes  \cdots \otimes
\varphi_{i_{2k}}) \right)I_{2l}  \left( (\varphi_{j_{1}}\otimes  \cdots \otimes
\varphi_{j_{2j}}) \right).
\end{eqnarray*}
The product formula (\ref{prod}) applied to the last equality yields
\begin{eqnarray*}
&& \sum_{i_{1}, \cdots , i_{2k}=1}^{n} a_{u, i_{1},\cdots , i_{2k}}\sum_{j_{1}, \cdots , j_{2l}=1}^{n} a_{u,j_{1},\cdots ,
j_{2l}}  I_{2k}  \left( (\varphi_{i_{1}}\otimes  \cdots \otimes
\varphi_{i_{2k}}) \right)I_{2l}  \left( (\varphi_{j_{1}}\otimes  \cdots \otimes
\varphi_{j_{2j}}) \right)\\
&=& \sum_{r=0}^{(2k)\wedge (2l)} r! C_{2k}^{r} C_{2l}^{r} \sum_{u_{1}, \cdots , u_{r}=1}^{n} \sum_{i_{1},\cdots , i_{2k-r} =1}^{n} \sum_{j_{1},\cdots , j_{2l-r} =1}^{n} a_{u, u_{1}, \cdots , u_{r}, i_{1}, \cdots , i_{2k-r}} a_{u, u_{1}, \cdots , u_{r}, j_{1}, \cdots , j_{2l-r}}\\
&&\times I_{2k+2l-2r} \left( \varphi_{i_{1}}\otimes \cdots \otimes \varphi _{i_{2k-r}}\otimes \varphi_{j_{1}}\otimes \cdots \otimes \varphi_{j_{2l-r}}\right)
\end{eqnarray*}
and therefore we obtain
\begin{eqnarray}
&& \langle DF_{n}, D(-L)^{-1}F_{n} \rangle \nonumber \\
&=& \sum_{k,l\geq 0} \frac{1}{(2k)! }\frac{1}{(2l+1)! }
\sum_{r=0}^{(2k)\wedge (2l)} r! C_{2k}^{r} C_{2l}^{r}\nonumber \\
&&\times  \sum_{u_{1}, \cdots , u_{r+1}=1}^{n} \sum_{i_{1},\cdots , i_{2k-r} =1}^{n} \sum_{j_{1},\cdots , j_{2l-r} =1}^{n}
 a_{u_{1}, u_{2}, \cdots , u_{r+1}, i_{1}, \cdots , i_{2k-r}} a_{u_{1}, u_{2}, \cdots , u_{r+1}, j_{1}, \cdots , j_{2l-r}}\nonumber \\
&&\times I_{2k+2l-2r} \left( \varphi_{i_{1}}\otimes \cdots \otimes \varphi _{i_{2k-r}}\otimes \varphi_{j_{1}}\otimes \cdots \otimes \varphi_{j_{2l-r}}\right).\label{end2}
\end{eqnarray}
\begin{remark}
\label{remarqueOrdreChaos}
The chaos of order zero in the above expression is obtained for $k=l$ and $r=2k$. It is therefore equal to
$$ \sum_{k\geq 0} \frac{1}{(2k)! }\frac{1}{(2k+1)! } (2k)! \sum_{i_{1},\cdots , i_{2k+1} =1}^{n} a_{i_{1},\cdots ,i_{2k+1} }^{2} $$
which is also equal to 1 as follows from Lemma \ref{lem:Chaosnorm}. Therefore it will vanish when we consider the difference $1-\langle DF_{n}, D(-L)^{-1}F_{n} \rangle $. This difference will have only chaoses of even orders.
\end{remark}

\noindent By changing the order of summation and by using the changes of indices $2k-r=k'$ and $2l-r=l'$, we can write
\begin{eqnarray*}
 &&\langle DF_{n}, D(-L)^{-1}F_{n} \rangle \\
 && =\sum_{r\geq 0}r! \sum_{2k\geq r} \sum_{2l\geq r} \frac{1}{(2k)! }\frac{1}{(2l+1)! }
 C_{2k}^{r} C_{2l}^{r}\\
 &&\times \sum_{u_{1}, \cdots , u_{r+1}=1}^{n} \sum_{i_{1},\cdots , i_{2k-r} =1}^{n} \sum_{j_{1},\cdots , j_{2l-r} =1}^{n}
 a_{u_{1}, u_{2}, \cdots , u_{r+1}, i_{1}, \cdots , i_{2k-r}} a_{u_{1}, u_{2}, \cdots , u_{r+1}, j_{1}, \cdots , j_{2l-r}}\\
&&\times I_{2k+2l-2r} \left( \varphi_{i_{1}}\otimes \cdots \otimes \varphi _{i_{2k-r}}\otimes \varphi_{j_{1}}\otimes \cdots \otimes \varphi_{j_{2l-r}}\right)\\
&=& \sum_{r\geq 0} \sum_{k,l\geq 0} \frac{1}{(k+r)! } \frac{1}{(l+r+1)!} C_{k+r}^{r} C_{l+r}^{r} \\
 &&\times \sum_{u_{1}, \cdots , u_{r+1}=1}^{n} \sum_{i_{1},\cdots , i_{k} =1}^{n} \sum_{j_{1},\cdots , j_{l} =1}^{n}
 a_{u_{1}, u_{2}, \cdots , u_{r+1}, i_{1}, \cdots , i_{k}} a_{u_{1}, u_{2}, \cdots , u_{r+1}, j_{1}, \cdots , j_{l}}\\
&&\times I_{2k+2l-2r} \left( \varphi_{i_{1}}\otimes \cdots \otimes \varphi _{i_{k}}\otimes \varphi_{j_{1}}\otimes \cdots \otimes \varphi_{j_{l}}\right)\\
&=&\sum_{k,l\geq 0}\sum_{r\geq 0} r!  \frac{1}{(k+r)! } \frac{1}{(l+r+1)!} C_{k+r}^{r} C_{l+r}^{r} \\
 &&\times \sum_{u_{1}, \cdots , u_{r+1}=1}^{n} \sum_{i_{1},\cdots , i_{k} =1}^{n} \sum_{j_{1},\cdots , j_{l} =1}^{n}
 a_{u_{1}, u_{2}, \cdots , u_{r+1}, i_{1}, \cdots , i_{k}} a_{u_{1}, u_{2}, \cdots , u_{r+1}, j_{1}, \cdots , j_{l}}\\
&&\times I_{k+l} \left( \varphi_{i_{1}}\otimes \cdots \otimes \varphi _{i_{k}}\otimes \varphi_{j_{1}}\otimes \cdots \otimes \varphi_{j_{l}}\right).
 \end{eqnarray*}
\noindent Once again using a change of indices ($k+l=m$), we obtain
\begin{eqnarray*}
 &&\langle DF_{n}, D(-L)^{-1}F_{n} \rangle \\
 &=&\sum_{m\geq 0} \sum_{k= 0}^{m} \sum_{r\geq 0}r!  \frac{1}{(k+r)! } \frac{1}{(m-k+r+1)!} C_{k+r}^{r} C_{m-k+r}^{r}\nonumber \\
 &&\times \sum_{u_{1}, \cdots , u_{r+1}=1}^{n} \sum_{i_{1},\cdots , i_{k} =1}^{n} \sum_{j_{1},\cdots , j_{m-k} =1}^{n}
 a_{u_{1}, u_{2}, \cdots , u_{r+1}, i_{1}, \cdots , i_{k}} a_{u_{1}, u_{2}, \cdots , u_{r+1}, j_{1}, \cdots , j_{m-k}}\nonumber \\
&&\times I_{m} \left( \varphi_{i_{1}}\otimes \cdots \otimes \varphi _{i_{k}}\otimes \varphi_{j_{1}}\otimes \cdots \otimes \varphi_{j_{m-k}}\right)\nonumber \\
&=&\sum_{m\geq 0} \sum_{k= 0}^{m} \frac{1}{k!} \frac{1}{(m-k)! } \sum_{r\geq 0} \frac{1}{r!} \frac{1}{m-k+r+1}\sum_{u_{1}, \cdots , u_{r+1}=1}^{n} \sum_{i_{1},\cdots , i_{m} =1}^{n} \nonumber \\
&&\times  a_{u_{1}, u_{2}, \cdots , u_{r+1}, i_{1}, \cdots , i_{k}} a_{u_{1}, u_{2}, \cdots , u_{r+1}, i_{k+1}, \cdots , i_{m}}I_{m} \left( \varphi_{i_{1}}\otimes \cdots \otimes \varphi _{i_{k}}\otimes \varphi_{i_{k+1}}\otimes \cdots \otimes \varphi_{i_{m}}\right)
 \end{eqnarray*}
where at the end we renamed the indices $j_{1},\cdots , j_{m-m}$ as $i_{k+1},\cdots , i_{m}$. We obtain
\begin{equation*}
\langle DF_{n}, D(-L)^{-1}F_{n} \rangle= \sum_{m\geq 0} I_{m} ( h_{m}^{(n)})
\end{equation*}
where
\begin{eqnarray}
h_{m}^{(n)} &=&\sum_{k= 0}^{m} \frac{1}{k!} \frac{1}{(m-k)! } \sum_{r\geq 0} \frac{1}{r!} \frac{1}{m-k+r+1}\sum_{u_{1}, \cdots , u_{r+1}=1}^{n} \sum_{i_{1},\cdots , i_{m} =1}^{n}\nonumber \\
&& a_{u_{1}, u_{2}, \cdots , u_{r+1}, i_{1}, \cdots , i_{k}} a_{u_{1}, u_{2}, \cdots , u_{r+1}, i_{k+1}, \cdots , i_{m}}\nonumber \\
 &&\varphi_{i_{1}}\otimes \cdots \otimes \varphi _{i_{k}}\otimes \varphi_{i_{k+1}}\otimes \cdots \otimes \varphi_{i_{m}}\label{hnm}
\end{eqnarray}
Let us make some comments about this result before going any further. These remarks will simplify the expression that we have just obtained. As follows from Lemma \ref{lem:intbypart}, the coefficients $a_{i_{1},\cdots ,i_{k}}$ are zero if $k$ is even. Therefore, the numbers $r+1+k$ and $r+1+m-k$ must be odd. This implies that $m$ must be even and this is coherent with our previous observation (see Remark \ref{remarqueOrdreChaos}) that the chaos expansion of $\langle DF_{n}, D(-L)^{-1}F_{n} \rangle$ only contains chaoses of even orders. The second comment concerns the chaos of order zero. If $m=0$ then $k=0$ and we obtain
\begin{equation*}
h_{0}^{(n)}= \sum_{r\geq 0}\sum_{u_{1}, \cdots , u_{r+1}=1}^{n} \frac{1}{r!} \frac{1}{r+1} a^{2}_{u_{1},\cdots , u_{r+1}} =\sum_{r\geq 1} \frac{1}{r!} \sum_{u_{1}, \cdots , u_{r}=1}^{n}  a^{2}_{u_{1},\cdots , u_{r}}.
\end{equation*}
Thus, because the summand $ \sum_{r\geq 1} \frac{1}{r!} \sum_{u_{1}, \cdots , u_{r}=1}^{n}  a^{2}_{u_{1},\cdots , u_{r}}-1$ is zero by using Lemma \ref{lem:Chaosnorm},
\begin{eqnarray*}
\langle DF_{n}, D(-L)^{-1}F_{n} \rangle -1&=& \left( \sum_{r\geq 1} \frac{1}{r!} \sum_{u_{1}, \cdots , u_{r}=1}^{n}  a^{2}_{u_{1},\cdots , u_{r}}-1\right)+ \sum_{m\geq 1} I_{2m}(h_{2m}^{(n)})\\
&=& \sum_{m\geq 1} I_{2m}(h_{2m}^{(n)})
\end{eqnarray*}
with $h_{2m}^{(n)}$ given by (\ref{hnm}).
\\~\\
\noindent Using the isometry formula of multiple integrals in order to compute the $L^{2}$ norm of the above expression and noticing that the function $h_{2m}^{(n)}$ is symmetric, we find that
\begin{eqnarray*}
&& \mathbf{E}\left( \left( \langle DF_{n}, D(-L)^{-1}F_{n} \rangle -1\right) ^{2}\right) = \sum_{m\geq 1}(2m)!\langle h_{2m}^{(n)}, h_{2m}^{(n)}\rangle _{L^{2}([0,1]^{2m})}\\
&=&\sum_{m\geq 1} (2m)!\sum_{k,l=0}^{2m} \frac{1}{k!}\frac{1}{l!} \frac{1}{(2m-k)!} \frac{1}{(2m-l)!} \sum_{r,q\geq 0} \frac{1}{r!}\frac{1}{q!}\frac{1}{2m-k+r+1}\frac{1}{2m-l+q+1}\\
&&\times \sum_{u_{1}, \cdots , u_{r+1}=1}^{n}\sum_{v_{1}, \cdots , v_{q+1}=1}^{n}\sum_{i_{1},\cdots , i_{2m} =1}^{n}\\
&&\times a_{u_{1}, u_{2}, \cdots , u_{r+1}, i_{1}, \cdots , i_{k}} a_{u_{1}, u_{2}, \cdots , u_{r+1}, i_{k+1}, \cdots , i_{2m}}a_{v_{1}, v_{2}, \cdots , v_{q+1}, i_{1}, \cdots , i_{k}} a_{v_{1}, v_{2}, \cdots , v_{q+1}, i_{k+1}, \cdots , i_{2m}}\\
&=&\sum_{m\geq 1}(2m)! \sum_{i_{1},\cdots , i_{2m} =1}^{n}\left( \sum_{k=0}^{2m} \frac{1}{k!} \frac{1}{(2m-k)!}\sum_{r\geq 0} \frac{1}{r!}\frac{1}{2m-k+r+1}
\sum_{u_{1}, \cdots , u_{r+1}=1}^{n}\right. \\
&&\left. a_{u_{1}, u_{2}, \cdots , u_{r+1}, i_{1}, \cdots , i_{k}} a_{u_{1}, u_{2}, \cdots , u_{r+1}, i_{k+1}, \cdots , i_{2m}}\right) ^{2},
\end{eqnarray*} which is the desired result.\qed \vskip0.3cm

\noindent Before proving our main result, let us discuss a particular case as an exemple in order to better understand the general phenomenon. This is both useful and important in order to have a good overview of the functioning of a simple case. Assume that $k=0$ and $l=1$. The corresponding summand in (\ref{end2}) reduces to
\begin{equation*}
\frac{1}{3!} \sum_{u=1}^{n} a_{u}\sum_{j_{1},j_{2}=1}^{n} a_{u,j_{1},j_{2}}I_{2}\left( \varphi_{j_{1}}\otimes \varphi_{j_{2}}\right).
\end{equation*}
Its $L^{2}$-norm is
\begin{equation*}
\frac{1}{3} \sum_{j_{1},j_{2}=1}^{n} \left( \sum_{u=1}^{n} a_{u}a_{u,j_{1},j_{2}}\right) ^{2}=\frac{1}{3}\sum_{j_{1}=1}^{n}\left( \sum_{u=1}^{n}a_{u}a_{u,j_{1},j_{1}}\right) ^{2}
\end{equation*}
because $a_{u,j_{1},j_{2}}=0$ if $j_{1}\not =j_{2}$. Using (\ref{imp}), it reduces to a quantity equivalent to
$$\frac{1}{3}( na_{1}^{2} a_{1,1,1}^{2} + n((n-1) a_{1} a_{1,1,2} )^{2})$$
which, using (\ref{imp}) again, is of order
$$n \left( \frac{1}{\sqrt{n}}\right) ^{2}\left(\frac{1}{n^{\frac{3}{2}}}\right) ^{2}+ n \left( (n-1) \frac{1}{\sqrt{n}}\frac{1}{n^{\frac{3}{2}}}\right)^{2} \sim n^{-1}.$$ The following theorem, which gathers all of the previous results of the paper, is the general equivalent of the toy exemple presented above.
\begin{thm*}\label{t1}
For any integer $n\geq 2$,
\begin{equation*}
\mathbf{E}\left( \left( \langle DF_{n}, D(-L)^{-1}F_{n} \rangle -1\right) ^{2}\right) \leq \frac{c_{0}}{n}
  \end{equation*}
  with
  \begin{eqnarray}
  \label{c0}
  c_{0}&=&\sum_{m\geq 1}(2m)!  \left(  \sum_{k=0}^{2m} \frac{1}{2k!} \frac{1}{(2m-2k)!}\sum_{r\geq 0} \frac{1}{(2r)!}\frac{1}{2m-2k+2r+1}c(k,r,m)\right)^{2}\\
&&+  \left( \sum_{k=0}^{2m} \frac{1}{(2k+1)!} \frac{1}{(2m-2k-1)!}\sum_{r\geq 0} \frac{1}{(2r-1)!}\frac{1}{2m-2k+2r+1}c(k,r,m)\right) ^{2}\nonumber
  \end{eqnarray}
and where $c(k,r,m)$ is given by (\ref{ckr}).
\end{thm*}
\noindent {\bf Proof: }Observe that the integers $r+1+k$ and $r+1 +2m-k$ both have to be odd numbers (otherwise the coefficients $a_{u_{1}, u_{2}, \cdots , u_{r+1}, i_{1}, \cdots , i_{k}}$ and $ a_{u_{1}, u_{2}, \cdots , u_{r+1}, i_{k+1}, \cdots , i_{2m}}$ vanish). This implies two cases: either $r$ is even and $k$ is even or $r$ is odd and $k$ is odd. Thus, we can write
\begin{eqnarray}
&& \mathbf{E}\left( \left( \langle DF_{n}, D(-L)^{-1}F_{n} \rangle -1\right) ^{2}\right) \nonumber\\
&=&\sum_{m\geq 1}(2m)! \sum_{i_{1},\cdots , i_{2m} =1}^{n}\left( \sum_{k=0}^{2m} \frac{1}{2k!} \frac{1}{(2m-2k)!}\sum_{r\geq 0} \frac{1}{(2r)!}\frac{1}{2m-2k+2r+1}
\sum_{u_{1}, \cdots , u_{2r+1}=1}^{n}\right. \nonumber\\
&&\left. a_{u_{1}, u_{2}, \cdots , u_{2r+1}, i_{1}, \cdots , i_{2k}} a_{u_{1}, u_{2}, \cdots , u_{2r+1}, i_{2k+1}, \cdots , i_{2m}}\right) ^{2}\nonumber\\
&&+\sum_{m\geq 1}(2m)! \sum_{i_{1},\cdots , i_{2m} =1}^{n}\left( \sum_{k=0}^{2m} \frac{1}{(2k+1)!} \frac{1}{(2m-2k-1)!}\sum_{r\geq 0} \frac{1}{(2r-1)!}\frac{1}{2m-2k+2r+1}
\sum_{u_{1}, \cdots , u_{2r}=1}^{n}\right. \nonumber\\
&&\left. a_{u_{1}, u_{2}, \cdots , u_{2r}, i_{1}, \cdots , i_{2k+1}} a_{u_{1}, u_{2}, \cdots , u_{2r}, i_{2k+2}, \cdots , i_{2m}}\right) ^{2}.\label{f1}
\end{eqnarray}
Let us treat the first part of the sum (\ref{f1}). Assume that the number of common numbers occurring in the sets $\{ u_{1},\cdots , u_{2r+1}\}$ and $\{ i_{1},\cdots , i_{2k}\}$ is $x$ and and the number of common numbers occurring in the sets $\{ u_{1},\cdots , u_{2r+1}\}$ and $\{i_{2k+1},\cdots ,i_{2m-2k}\} $ is $y$. This can be formally written as
$$ \left| \{  u_{1},\cdots , u_{2r+1}\} \cap\{  \{ i_{1},\cdots , i_{2k}\}\right|  =x$$
and
$$ \left| \{ u_{1},\cdots , u_{2r+1}\} \cap\{  i_{2k+1},\cdots ,i_{2m-2k}\}\right|  =y.$$
It is clear that
$$x\leq (2r+1) \wedge 2k \mbox{ and } y\leq (2r+1) \wedge 2m-2k.$$
This also implies $x+y\leq 2m$. According to the definitions of $x$ and $y$, it can be observed that $x$ and $y$ must be even. We will denote them by $2x$ and $2y$ from now on.

\noindent The next step in the proof is to determine how many distinct sequences of numbers can occur in the set
$$\{ u_{1},\cdots , u_{2r+1}, i_{1}, \cdots ,i_{2k} \}.$$
We can have sequences of lengths (all of the lengths that we consider from now on are greater or equal to one)
 $2c_{1},2c_{2},\cdots , 2c_{l_{1}}$ with $2(c_{1}+\cdots +c_{l_{1}})=2x$ in the set $\{u_{1},\cdots ,u_{2r+1} \} \cap \{i_{1},\cdots ,i_{2k}\}$
but also sequences of lengths
  $2d_{1}, 2d_{2},\cdots , 2d_{l_{2}}$ with $2( d_{1}+\cdots +d_{l_{2}}) =2k-2x$ in the set $\{i_{1},\cdots ,i_{2k}\}\setminus \{u_{1},\cdots ,u_{2r+1} \}$
as well as sequences of lengths
  $2e_{1}+1, 2e_{2},\cdots , 2e_{l_{3}}$ with $1+2(e_{1}+\cdots +e_{l_{3}})= 2r+1-2x$ in the set $\{u_{1},\cdots ,u_{2r+1} \}\setminus \{i_{1},\cdots ,i_{2k}\}.$
In this last sequence we have one (and only one) length equal to 1 (because we are allowed to choose only one odd number in the set  $ \{u_{1},\cdots ,u_{2r+1} \}\setminus \{i_{1},\cdots ,i_{2k}\}$). We will have, if we have a configuration as above,
 \begin{equation*}
 a_{u_{1}, u_{2}, \cdots , u_{2r+1}, i_{1}, \cdots , i_{2k}}\leq c(r,c,e) n ^{-\frac{1}{2} -l_{1}-l_{2}-l_{3} }
 \end{equation*}
 where
 \begin{equation}
 \label{cl1l2}c(r,c,e)=r! (2r-1)!! \frac{ (2c_{1})!\cdots (2c_{l_{1}})!(2e_{1}+1)! (2e_{2})!\cdots (2e_{l_{3}})! }{(c_{1}! \cdots c_{l_{1}}! e_{1}! \cdots e_{l_{3}}!)^{2}}t(c_{1})\cdots t(c_{l_{1}})t(e_{1})\cdots t(e_{l_{3}})
 \end{equation}
 and the constants $t$ are given by (\ref{tdj}).

\noindent In the same way, assuming that we have sequences of lengths
 $2f_{1},2f_{2},\cdots , 2f_{l_{4}}$ with $2(f_{1}+\cdots +f_{l_{4}})= 2m-2k-2y$ in the set $\{i_{2k+1},\cdots , i_{2m} \} \setminus \{u_{1},\cdots ,u_{2r+1} \}$
 and sequences of lengths $2g_{1}+1, 2g_{2},\cdots , 2g_{l_{5}}$ with $1+2(g_{1}+\cdots +g_{5})=2r+1-2y$ in the set $\{u_{1},\cdots ,u_{2r+1}\} \setminus \{i_{2k+1},\cdots ,i_{2m}\}.$
We will obtain
 \begin{equation*}
 a_{u_{1}, u_{2}, \cdots , u_{2r+1}, i_{2k+1}, \cdots , i_{2n}}\leq c(k, c,d)n ^{-\frac{1}{2} -l_{1}-l_{4}-l_{5} +1}
 \end{equation*}
 with $c(k, c,d)$  defined as in (\ref{cl1l2}).
 The sum over $u_{1},\cdots ,u_{r+1} $ from $1$ to $n$ reduces to a sum of $l_{1}+l_{3}+l_{5}-1 $ distinct indices from $1$ to $n$. Therefore we get
\begin{eqnarray*}
&&\sum_{u_{1}, \cdots , u_{2r+1}=1}^{n}
 a_{u_{1}, u_{2}, \cdots , u_{2r+1}, i_{1}, \cdots , i_{2k}} a_{u_{1}, u_{2}, \cdots , u_{2r+1}, i_{2k+1}, \cdots , i_{2n}}\\
 &&\leq c(k,r,m)n^{-l_{1}- l_{2}-l_{4}}
\end{eqnarray*}
with
\begin{equation}\label{ckr}
c(k,r,m)=\sum_{x+y=2m}\sum_{c_{1}+\cdots +c_{l_{1}}=x}\sum_{d_{1}+\cdots +d_{l_{2}}=y}\sum_{e_{1}+\cdots +e_{l_{3}}=r-x}c(r,c,e)c(k, c, d).
\end{equation}
We need to consider the sum $i_{1}, \cdots , i_{2m}$ from $1$ to $n$. It reduces to a sum over $l_{2}+l_{4}$ distinct indices. Thus
\begin{eqnarray*}
&& \sum_{i_{1},\cdots , i_{2m} =1}^{n}\left( \sum_{k=0}^{2m} \frac{1}{2k!} \frac{1}{(2m-2k)!}\sum_{r\geq 0} \frac{1}{(2r)!}\frac{1}{2m-2k+2r+1}
\sum_{u_{1}, \cdots , u_{2r+1}=1}^{n}\sum_{u_{1}, \cdots , u_{2r+1}=1}^{n}\right. \\
&&\left. a_{u_{1}, u_{2}, \cdots , u_{2r+1}, i_{1}, \cdots , i_{2k}} a_{u_{1}, u_{2}, \cdots , u_{2r+1}, i_{2k+1}, \cdots , i_{2m}}\right) ^{2}\\
&&\leq n^{l_{2}+l_{4}}\left( \frac{1}{n^{2l_{1}+l_{2} + l_{4} }}\right) ^{2}  \sum_{k=0}^{2m} \frac{1}{2k!} \frac{1}{(2m-2k)!}\sum_{r\geq 0} \frac{1}{(2r)!}\frac{1}{2m-2k+2r+1}c(k,r,m)\\
&=&\frac{1}{n^{2l_{1}+l_{2}+l_{4}}}  \left(  \sum_{k=0}^{2m} \frac{1}{2k!} \frac{1}{(2m-2k)!}\sum_{r\geq 0} \frac{1}{(2r)!}\frac{1}{2m-2k+2r+1}c(k,r,m)\right)^{2}.
\end{eqnarray*}

\noindent Note that either $l_{1}+l_{2}\geq 1$ or $l_{1}+l_{4} \geq 1$ (this is true because $m\geq 1$).  Then this term is at most  of order of $n^{-1}$.

\noindent Let us now look at the second part of the sum in (\ref{f1}). Suppose that in the sets
$\{u_{1},\cdots , u_{2r}\} \cap \{i_{1},\cdots ,i_{2k+1}\}$, $\{i_{1},\cdots ,i_{2k+1}\}\setminus \{u_{1},\cdots , u_{2r}\}$, $\{u_{1},\cdots , u_{2r}\}\setminus \{i_{1},\cdots ,i_{2k+1}\}$, $\{ i_{2k+2},\cdots , i_{2m-2k}\}\setminus \{u_{1},\cdots , u_{2r}\}$, $\{u_{1},\cdots , u_{2r}\}\setminus \{ i_{2k+2},\cdots , i_{2m-2k}\}$
we have sequences with lengths
$$p_{1}, \hskip0.2cm p_{2}, \hskip0.2cm p_{3}, \hskip0.2cm p_{4}, \hskip0.2cm p_{5} \hskip0.2cm$$
respectively (the analogous of $l_{1},\cdots , l_{5}$ above). Then the behavior with respect to $n$ of
$$\sum_{u_{1}, \cdots , u_{2r}=1}^{n}
 a_{u_{1}, u_{2}, \cdots , u_{2r}, i_{1}, \cdots , i_{2k+1}} a_{u_{1}, u_{2}, \cdots , u_{2r}, i_{2k+2}, \cdots , i_{2m}}$$
 is of order of $n^{p_{1}+p_{3}}\frac{1}{ n^{2p_{1} + p_{3} +p_{4} }}$. Therefore the behavior with respect to $n $ of the second sum in (\ref{f1}) is of order
 $$n^{p_{2}+1+p_{4}+ 1} \left( \frac{1}{n^{1}+2p_{1}+p_{2}+p_{4}} \right)^{2}= \frac{1}{n^{2p_{1}+p_{2}+p_{4}}}.$$
 Again, since either $p_{1}+p_{2}\geq 1$ or $p_{1}+p_{4} \geq 1$, the behavior of the term is at most of order $n^{-1}$. Therefore
\begin{equation*}
\mathbf{E}\left( \left( \langle DF_{n}, D(-L)^{-1}F_{n} \rangle -1\right) ^{2}\right) \leq \frac{c_{0}}{n}
\end{equation*}
where the constant $c_{0}$ is given by (\ref{c0}). The fact that the sum over $m $ is finite is a consequence of the following argument: $\langle DF_{n}, D(-L)^{-1}F_{n} \rangle$ belongs to $\mathbb{D}^{\infty , 2}(\Omega)$ (which is true  based on the derivation rule - Exercise 1.2.13 in \cite{N}- and since $F_{n}$ belongs to $\mathbb{D} ^{\infty, 2}$ as a consequence of Proposition 1.2.3 in \cite{N}), this implies that $\sum_{m} m! m^{k} \Vert h_{m}^{(n)} \Vert _{2}^{2} <\infty$ for every $k$ where $h_{m}^{(n)} $ is given by (\ref{hnm}). Therefore, the constant $c(m,k, r)$ defined in (\ref{ckr}) behaves at most as a power function with respect to $m.$  \qed

\vskip0.3cm
\begin{corollary}
Let $J_{m}(F_{n})$ denotes the projection on the $m^{\mbox{\tiny{th}}}$ Wiener chaos of the random variable $F_{n}$. Then for every $m\geq 1$ the sequence $J_{m}(F_{n})$ converges as $n\to \infty$ to a standard normal random variable.
\end{corollary}
\noindent {\bf Proof: }The proof is a consequence of the proof of Theorem \ref{t1}. \qed
\vskip0.5cm

{\bf Acknowledgement: } The authors wish to thank Natesh Pillai for interesting discussions and for his contribution to some parts of this work. The second  author is partially supported by the ANR grant "Masterie" BLAN 012103.

\end{document}